\newcommand{\pbaddress}{yochay.jerby@gmail.com}

\documentclass[12pt]{amsart}

\usepackage{epsfig}
\usepackage{subfig}
\usepackage{amscd}
\usepackage[mathscr]{eucal}
\usepackage{amssymb}
\usepackage{amsxtra}
\usepackage{amsmath}
\usepackage[all]{xy}
\usepackage{tikz}
\usepackage{pgf}
\usetikzlibrary{calc,backgrounds,arrows,matrix,shapes}
\usepackage{verbatim}
\usepackage[bookmarks]{hyperref}

\theoremstyle{plain}

\theoremstyle{definition}

\theoremstyle{remark}

\renewcommand{\themainthm}{\Alph{mainthm}}
\newcommand{\cntrs}{\setcounter{thm}{0}
  \renewcommand{\thethm}{\thesection.\Alph{thm}}}
\newcommand{\cntrsb}{\setcounter{thm}{0}
  \renewcommand{\thethm}{\thesubsection.\Alph{thm}}}

\newcommand{\Qed}{\hfill \qedsymbol \medskip}
\newcommand{\INS}{{\noindent{\large\textbf{INSERT HERE: }}}}
\newcommand{\INSREF}{{\large{\textbf{INSERT REF.}}}}
\newcommand{\NOTE}{\noindent{\large\textbf{NOTE !!!}}}

\newcommand{\Id}{{{\mathchoice {\rm 1\mskip-4mu l} {\rm 1\mskip-4mu l}
      {\rm 1\mskip-4.5mu l} {\rm 1\mskip-5mu l}}}}

\renewcommand{\baselinestretch}{1.2}

\begin{document}

\title{On Landau-Ginzburg systems, Quivers and Monodromy}

\address{Institut de Mathematiques, Rue Emile-Argand 11 CH-2000 Neuch$\hat{\textrm{a}}$tel}
\email{yochay.jerby@unine.ch}

\date{\today}

\author{Yochay Jerby}

%\bibliographystyle{plain}

%----------------------------------------------------------------------
%
% Abstract
%
\begin{abstract}
Let $X$ be a toric Fano manifold and denote by $Crit(f_X) \subset (\mathbb{C}^{\ast})^n$ the solution scheme of the corresponding Landau-Ginzburg system of equations. For toric Del-Pezzo surfaces and various toric Fano threefolds we define a map $L : Crit(f_X) \rightarrow Pic(X)$ such that $\mathcal{E}_L(X) : = L(Crit(f_X)) \subset Pic(X)$ is a full strongly exceptional collection of line bundles.  We observe the existence of a natural monodromy map $$ M : \pi_1(L(X) \setminus R_X,f_X) \rightarrow Aut(Crit(f_X))$$ where $L(X)$ is the space of
all Laurent polynomials whose Newton polytope is equal to the Newton polytope of $f_X$, the Landau-Ginzburg potential of $X$, and $R_X \subset L(X)$ is the space of all elements whose corresponding solution scheme is reduced. We show that monodromies of $Crit(f_X)$
admit non-trivial relations to quiver representations of the exceptional collection $\mathcal{E}_L(X)$. We refer to this property as the $M$-aligned property of the maps $L: Crit(f_X) \rightarrow Pic(X)$.  We discuss possible applications of the existence of such $M$-aligned exceptional maps
to various aspects of mirror symmetry of toric Fano manifolds.

\end{abstract}

\maketitle

%----------------------------------------------------------------------
%
% Beginning of text
%

\section{Introduction and Summary of Main Results}
\label{s:intro}

\hspace{-0.6cm} Let $X$ be a smooth algebraic manifold and let $\mathcal{D}^b(X)$ be the bounded derived category of coherent sheaves on $X$, see \cite{GM}. Let $A$ be a finite dimensional associative algebra over the complex numbers and
let $\mathcal{D}^b(A)$ be the derived category of right modules over $A$. One of the fundamental questions in the study of $\mathcal{D}^b(X)$ is the following: \emph{Given a manifold $X$, is $\mathcal{D}^b(X)$ equivalent to the derived category $\mathcal{D}^b(A)$ of some finite dimensional associative algebra $A$?}

\hspace{-0.6cm} The first example of such an equivalence is Beillinson's famous description of $\mathcal{D}^b(X)$ for $X=\mathbb{P}^n$, see \cite{B}. Beilinson shows that $\mathcal{D}^b(\mathbb{P}^n)$ is equivalent to $\mathcal{D}^b(A_n)$ where
$A_n=End(T_n)$ is the endomorphism ring of the vector bundle $$ T_n = \mathcal{O} \oplus \mathcal{O}(1) \oplus ... \oplus \mathcal{O}(n)$$ In general, an object $ E \in \mathcal{D}^b(X)$ is said to be \emph{exceptional} if $Hom(E,E)=\mathbb{C}$ and $Ext^i(E,E)=0$ for $0<i$. An ordered collection $\left \{ E_1,...,E_N \right \} \subset \mathcal{D}^b(X)$
is said to be an \emph{exceptional collection} if each $E_j$ is exceptional and $$Ext^i(E_k,E_j) =0 \textrm{ for } j<k \textrm{ and } 0 \leq i $$ An
 exceptional collection is said to be \emph{strongly exceptional} if also $Ext^i(E_j,E_k)=0$ for $j \leq k$ and $0<i$. A strongly exceptional collection is called \emph{full} if its elements generate $\mathcal{D}^b(X)$ as a triangulated category. In particular, if $ \mathcal{E} = \left \{ E_1,...,E_N \right \} \subset \mathcal{D}^b(X)$
is a full strongly exceptional collection of objects, the corresponding adjoint functors $$\begin{array}{ccc} R Hom_X(T, -) : \mathcal{D}^b(X) \rightarrow \mathcal{D}^b(A_T) & ; & - \otimes^L_{A_T} T : \mathcal{D}^b(A_T) \rightarrow \mathcal{D}^b(X) \end{array} $$ are equivalences of categories,
where $T= \bigoplus_{i=1}^N E_i$. For a given algebraic manifold $X$ one thus asks the following two, related, but not similar questions: (a) does $X$ admit a full exceptional collection of objects in $\mathcal{D}^b(X)$? (b) does $X$ admit a full strongly exceptional collection of \emph{line bundles} in $Pic(X)$?

\hspace{-0.6cm} A class of manifolds on which these questions have been extensively studied in recent years is the class of toric manifolds and, specifically, the class
of toric Fano manifolds, see \cite{BH,BT,CMR,CMR2,CMR3,Ka,K,Pe,U}. Question (a) was answered affirmatively by Kawamata which showed that \emph{any toric manifold admits a full exceptional collection of objects in $\mathcal{D}^b(X)$}, see \cite{Ka}. However, the
more refined question (b) of which toric manifolds admit full strongly exceptional collections of line bundles in $Pic(X)$ is currently completely
open.

\hspace{-0.6cm} Question (b) has been especially studied for the class of toric Fano manifolds and, indeed, many examples of toric Fano manifolds which admit full strongly exceptional collections have been discovered by various authors.
The abundance of examples led experts to ask whether, in fact, any toric Fano manifold admits a full strongly exceptional
collection of line bundles in $Pic(X)$, see \cite{BH,CDRMR}. However, in a recent surprising work \cite{E} Efimov discovered examples of
toric Fano manifolds which do not admit any full strongly exceptional collections of line bundles. In particular, the question of which toric Fano manifolds admit full strongly exceptional collections in $Pic(X)$ is currently still open.

\hspace{-0.6cm} On the other hand, the theory of quantum cohomology introduces a family of commutative associative operations $\circ_{\omega} : H^{\ast}(X) \otimes H^{\ast}(X) \rightarrow H^{\ast}(X)$ parameterized by
classes $\omega \in H^{\ast}(X)$. This family of "quantum products" defines the structure of a Frobenius super-manifold over $H^{\ast}(X)$, which is known as the \emph{big quantum cohomology} of $X$. In particular, the big quantum cohomology is said to be
\emph{semi-simple} if the operation $\circ_{\omega}$ is a semi-simple ring operation for generic $\omega \in H^{\ast}(X)$. One of the fundamental conjectures on the structure of $\mathcal{D}^b(X)$ is the Dubrovin-Bayer-Manin conjecture, which relates the existence of full  exceptional collections of objects in $\mathcal{D}^b(X)$ to the semi-simplicty of the
big quantum cohomology of $X$, see \cite{BM,Du}.

\hspace{-0.6cm} When $X$ is a toric manifold, the Dubrovin-Bayer-Manin conjecture is actually known to hold due to the combined results of Kawamata (on the existence of full exceptional collections of
objects in $\mathcal{D}^b(X)$, see \cite{Ka}) and Iritani (on the semi-simplicity of the big quantum cohomology of toric manifolds, see \cite{I}). In view of the above one is led to ask whether it is possible to relate further, more refined, properties of quantum cohomology to the existence of full strongly exceptional collections of line bundles in $Pic(X)$?

\hspace{-0.6cm} Indeed, of special importance in quantum cohomology theory is the fiber $QH(X) \simeq (H^{\ast}(X),\circ_0)$ of the
 big quantum cohomology over $\omega=0$, which is known as the \emph{small}
 quantum cohomology ring of $X$. When $X$ is a toric Fano manifold the small quantum cohomology is expressed as the Jacobian
 ring of the Landau-Ginzburg potential, which is a Laurent polynomial $f_X \in \mathbb{C}[z_1^{\pm},...,z_n^{\pm}]$ associated to $X$, see \cite{Ba,FOOO,OT}. Consider the system of algebraic equations $ z_i \frac{\partial}{\partial z_i } f_X(z_1,...,z_n)=0$ for $i=1,...,n $, to which we refer as the Landau-Ginzburg system of equations of $X$ and denote by $Crit(f_X) \subset (\mathbb{C}^{\ast})^n$ the corresponding solution scheme. Our aim in this work is to present, via examples, various relations between properties of the solution scheme $Crit(f_X) \subset (\mathbb{C}^{\ast})^n$ and properties of full strongly exceptional collections of line bundles $\mathcal{E} \subset Pic(X)$. In particular, these relations lead us to suggest a "small variation" of the Dubrovin-Bayer-Manin conjecture for toric Fano manifolds, which we formulate below. As a starting point, consider the following example:

\bigskip

\hspace{-0.6cm} \bf Example \rm (projective space): For $X= \mathbb{P}^n$ the Landau-Ginzburg potential is given by $f(z_1,...,z_n)=z_1+...+z_n +\frac{1}{z_1 \cdot ... \cdot z_n}$ and the corresponding system
of equations is $$ z_i \frac{\partial}{\partial z_i } f_X(z_1,...,z_n)=z_i - \frac{1}{z_1 \cdot ... \cdot z_n} =0 \hspace{0.5cm} \textrm{ for } \hspace{0.25cm} i=1,...,n $$ The solution scheme $Crit(f_X) \subset (\mathbb{C}^{\ast})^n$ is given by
$z_k= ( e^{\frac{2 \pi ki}{n+1}},...,e^{\frac{2 \pi k i}{n+1}}) $ for $k=0,...,n$.

\bigskip

\hspace{-0.6cm} In general, Ostrover and Tyomkin show in \cite{OT} that $X$ has semi-simple
quantum cohomology if and only if the number of elements of $Crit(f_X)$ is $\chi(X)$, the Euler characteristic of $X$. On the other hand, the expected number of elements in a full strongly exceptional collection in $Pic(X)$ is also $\chi(X)$, see \cite{E}. In view of this, we refer to a map $L: Crit(f_X) \rightarrow Pic(X)$ as an \emph{exceptional map} if its image $\mathcal{E}_L(X) = L(Crit(f_X)) \subset Pic(X)$ is a full strongly exceptional collection.
The guiding question is thus the following:

\bigskip

\hspace{-0.6cm} \bf Main Question \rm (small toric Fano DBM-conjecture): Does any toric Fano manifold $X$ whose small quantum cohomology
$QH(X)$ is semi-simple admit an exceptional map $L : Crit(f_X) \rightarrow Pic(X)$ naturally generalizing the map $L(z_k)=\mathcal{O}(k)$ in the
case of projective space?

\bigskip

\hspace{-0.6cm}  Note that defining an exceptional map $L : Crit(f_X) \rightarrow Pic(X)$ requires the association of integral invariants to elements of $Crit(f_X)$. In the case of projective space, such an association is, in fact, misleadingly simple as the entries of
the elements are given
as roots of unity. However, in general, this is not the case, which is one of the main difficulties in defining the exceptional maps in full generality. Instead, in this work, we consider a few specific examples of toric Fano manifolds, which are known to admit full strongly
exceptional collections of line bundles. The manifolds considered are the following: (a) toric Del-Pezzo surfaces, (b) Fano $\mathbb{P}^1$-bundles over $\mathbb{P}^2$, (c) Fano $\mathbb{P}^2$-bundles over $\mathbb{P}^1$, (d) $\mathbb{P}^1$-bundles over $\mathbb{P}^1 \times \mathbb{P}^1$.

\hspace{-0.6cm} Based on a study of the properties of $Crit(f_X)$, we show that the manifolds (a)-(d) could be naturally introduced with exceptional maps $L : Crit(f_X) \rightarrow Pic(X)$, generalizing the above
example of projective space. Our main result, however, is that once the mentioned exceptional maps are defined, various algebraic properties of the corresponding collection $\mathcal{E}_L(X)$ turn to be non-trivially related to geometric properties of
the solution scheme $Crit(f_X)$.

\hspace{-0.6cm} Indeed, one of the fundamental features of exceptional collections is their relation to quiver representations. Recall that a \emph{quiver with relations} $\widetilde{Q}=(Q,R)$ is a directed graph $Q$ with a two sided ideal $R$ in the path algebra $\mathbb{C}Q$ of $Q$, see \cite{DW}. In particular, a quiver with relations $\widetilde{Q}$ determines the following
associative algebra $A_{\widetilde{Q}}=\mathbb{C}Q/R$, called the path algebra of $\widetilde{Q}$. A standard construction associates to a collection of elements $\mathcal{C} \subset \mathcal{D}^b(X)$ and a basis $ B \subset A_{\mathcal{C}}$
a quiver with relations
$\widetilde{Q}(\mathcal{C},B)$ whose vertex set is $\mathcal{C}$ such that $A_{\mathcal{C}} \simeq A_{\mathcal{Q}(\mathcal{C},B)}$, see \cite{K}. 

   \hspace{-0.6cm} In our case, we observe that the algebras $A_{\mathcal{E}_L(X)}$ admit a natural basis $B_{\mathcal{E}_L(X)}$, which is 
uniquely
determined by the toric data. To the toric Fano manifolds (a)-(d) we thus associate the quiver with relations given by $\widetilde{Q}(\mathcal{E}_L(X)):=\widetilde{Q}(\mathcal{E}_L(X),B_{\mathcal{E}_L(X)})$. In the toric Del-Pezzo case, that is case (a), the quivers are similar to the ones described in \cite{K,Pe}. Moreover, due to the construction, the edge set of the quivers, $\widetilde{Q}_1(\mathcal{E}_L(X))$, is endowed with a map of the form $D : \widetilde{Q}_1(\mathcal{E}_L(X)) \rightarrow Div_T(X)$, where $Div_T(X)$ is the space of toric divisors of $X$ (see section 3). For a divisor $D \in Div_T(X)$ denote by $ \widetilde{Q}^D(\mathcal{E}_L(X))$ the sub-quiver of $ \widetilde{Q}(\mathcal{E}_L(X))$
whose edges satisfy $D(a)=D$.

\hspace{-0.6cm} On the other hand, in the quantum cohomology side, we note that there exists a natural mondromy group action on the solution set $Crit(f_X) \subset (\mathbb{C}^{\ast})^n$. Indeed, let $L(X) \subset \mathbb{C}[z_1^{\pm},...,z_n^{\pm}]$ be the vector space of Laurent polynomials whose Newton polytope is the same as that of $f_X$. As for $f_X$,
one can associate a scheme $Crit(f)$ to any element $f \in L(X)$. For a generic $f \in L(X)$ the scheme $Crit(f)$ is given as
the solution scheme of the system of equations $ z_i \frac{\partial}{\partial z_i } f(z_1,...,z_n)=0$ for $i=1,...,n $. Consider the hypersurface $R_X \subset L(\Delta^{\circ})$ of all $f \in L(X)$ such that $Crit(f)$ is non-reduced. Hence, in particular, when $QH(X)$ is semisimple one obtains, via standard analytic continuation, a monodromy map of the following form $$ M : \pi_1(L(X)
\setminus R_X, f_X) \rightarrow Aut(Crit(f_X))$$ For any manifold $X$ of (a)-(d) we describe a map $ \Gamma: Div_T(X) \rightarrow \pi_1(L(X) \setminus R_X , f_X) $, for the definition of this map see section 5. Denote by $\widetilde{Q}(D)$ the quiver whose vertex set is $\widetilde{Q}_0(D) =Crit(f_X)$, and which has an edge $a_{D}(z) \in
\widetilde{Q}_1(D)$ beginning at $z$ and ending at $M(\Gamma(D))(z)$, for any $z \in Crit(f_X)$.  Our main result is the following:

\bigskip

\hspace{-0.6cm} \bf Theorem A \rm ($M$-aligned property): Let $X$ be a toric Fano manifold of (a)-(d) and let $L: Crit(f_X) \rightarrow Pic(X)$ be the corresponding exceptional map. Then $$ \widetilde{Q}^{D}(\mathcal{E}_L(X)) \subset \widetilde{Q}(D) \hspace{0.5cm}
 \textrm{for any} \hspace{0.25cm} D \in Div_T(X)$$

\hspace{-0.6cm} We refer to the property described in Theorem A as the $M$-aligned property. Let us conclude by noting that although, by definition, any toric Fano manifold which admits a full strongly exceptional collection can be endowed with an exceptional map, the
 existence of $M$-aligned exceptional maps is less trivial. In particular, we suggest that the existence of $M$-aligned maps indicates that the relations between elements of $Crit(f_X)$ and full strongly exceptional collections of line bundles, described for manifolds (a)-(d), could be generalized to other examples of toric Fano manifolds with semi-simple quantum cohomology.

\bigskip

\hspace{-0.6cm} The rest of the work is organized as follows: In section 2 we recall relevant facts on toric Fano manifolds. In section 3 we review relevant aspects of the theory of derived categories and exceptional collections,
give examples of full strongly exceptional collections of line bundles and present their corresponding quivers. In section 4
we introduce the monodromy operator $M$. In section 5 we define the exceptional maps and compute their corresponding monodromies. Section 6 is devoted for concluding remarks and discussion of further relations to mirror symmetry.

\section{Relevant Facts on Toric Fano Manifolds}
\label{s:Rfotfm}

\hspace{-0.6cm} In this section we review relevant facts on toric Fano manifolds, we refer the reader to \cite{F,O} for a detailed overview of the theory of toric geometry.
A toric variety is an algebraic variety $X$ containing an algebraic torus $T \simeq (\mathbb{C}^{\ast})^n$ as a dense subset such that the action of $T$ on itself extends to the whole variety. A compact toric variety $X$ is said to be Fano if its anticanonical class $-K_X$ is Cartier and ample. In \cite{Ba2} Batyrev shows that there is a one to one correspondence between toric Fano varieties and reflexive polytopes.

\hspace{-0.6cm} Let $N \simeq \mathbb{Z}^n$ be a lattice and let $M = N^{\vee}=Hom(N, \mathbb{Z})$ be the dual lattice. Denote by $N_{\mathbb{R}} = N \otimes \mathbb{R}$ and $M_{\mathbb{R}}=M \otimes \mathbb{R}$ the corresponding vector space. Let $ \Delta \subset M_{\mathbb{R}}$ be an integral polytope and let $$ \Delta^{\circ} = \left \{ n \mid (m,n) \geq -1 \textrm{ for every } m \in \Delta \right \} \subset N_{\mathbb{R}}$$ be the \emph{polar} polytope of $\Delta$.

\bigskip

\hspace{-0.6cm} \bf Definition 2.1: \rm The polytope $\Delta \subset M_{\mathbb{R}}$ is said to be \emph{reflexive} if $0 \in \Delta$ and $\Delta^{\circ} \subset N_{\mathbb{R}}$ is integral. A reflexive polytope $\Delta$ is said to be \emph{Fano} if every facet of $\Delta^{\circ}$ is the convex hall of a basis of $M$.

\bigskip

\hspace{-0.6cm} In general, let $\Delta \subset M_{\mathbb{R}}$ be a polytope and let $$ L(\Delta) = \bigoplus_{m \in \Delta \cap M} \mathbb{C} m $$ be the space of Laurent polynomials whose Newton polytope is $\Delta$. A polytope $\Delta \subset M_{\mathbb{R}} $ determines an embedding $i_{\Delta} : (\mathbb{C}^{\ast})^n \rightarrow \mathbb{P}(L(\Delta)^{\vee})$ given by $ z \mapsto [z^m \mid m \in \Delta \cap M] $. The \emph{toric variety} $X_{\Delta} \subset \mathbb{P}(L(\Delta)^{\vee})$ corresponding to the polytope $\Delta \subset M_{\mathbb{R}}$ is defined
to be the compactification of the embedded torus $i_{\Delta}((\mathbb{C}^{\ast})^n) \subset \mathbb{P}(L(\Delta)^{\vee})$. It is shown by Batyrev in \cite{Ba2} that $X_{\Delta}$ is a Fano variety if $\Delta$ is reflexive and, in this case, the embedding
$i_{\Delta}$ is the anti-canonical embedding. The Fano variety $X_{\Delta}$ is smooth if and only if $\Delta^{\circ}$ is a Fano polytope.

\hspace{-0.6cm} Batyrev shows in \cite{Ba3} that there are a finite number of reflexive polytopes of given dimension. In particualr, there are a finite number of Fano polytopes in given dimension. In dimension three there are, up to equivalence, 18 Fano polytopes, see \cite{Ba4,Wa}. There are 124 Four dimensional Fano polytopes, up to equivalence. 123 of them were classified by Batyrev in \cite{Ba3} and an additional example was discovered by Sato in \cite{Sa}. The five dimensional Fano polytopes (866 examples, up to equivalence) were recently classified by Kreuzer and Nill in \cite{KN}. The two dimensional, Del-Pezzo, case is described in
the following example:

\bigskip

\hspace{-0.6cm} \bf Example 2.2 \rm (Toric Del-Pezzo surfaces): Up to integral automorphisms of $\mathbb{R}^2$ there are five Fano polytopes in dimension two, also called Del-Pezzo polytopes. Denote by $N_{\Delta^{\circ}}$ the matrix of vertices of $\Delta^{\circ}$.

\bigskip

\hspace{-0.6cm} (1) $ N_{\Delta^{\circ}} = \left ( \begin{array}{ccc} 1 & 0 & -1 \\ 0 & 1 & -1 \end{array} \right ) $ the corresponding manifold is $X_{\Delta}=\mathbb{P}^2$ .

\bigskip

\hspace{-0.6cm} (2) $ N_{\Delta^{\circ}}=\left ( \begin{array}{cccc} 1 & 0 & -1 & 0 \\ 0 & 1 & 0 & -1 \end{array} \right ) $ the corresponding manifold is $X_{\Delta}=\mathbb{P}^1 \times \mathbb{P}^1$.

\bigskip

\hspace{-0.6cm} (3) $ N_{\Delta^{\circ}}=\left ( \begin{array}{cccc} 1 & 0 & -1 & -1 \\ 0 & 1 & 0 &-1 \end{array} \right ) $ the corresponding manifold is $X_{\Delta}=Bl_1(\mathbb{P}^2)$, the blow up of projective plane
at one $T$-equivariant point.

\bigskip

\hspace{-0.6cm} (4) $ N_{\Delta^{\circ}}=\left ( \begin{array}{ccccc} 1 & 0 & -1 & 0 & -1 \\ 0 & 1 & 0 & -1 & -1 \end{array} \right ) $ the corresponding manifold is $X_{\Delta}=Bl_2(\mathbb{P}^2)$, the blow up of projective plane
at two $T$-equivariant points.

\bigskip

\hspace{-0.6cm} (5) $N_{\Delta^{\circ}}=\left ( \begin{array}{cccccc} 1 & 0 & 1 & -1 & 0 & -1 \\ 0 & 1 & 1 & 0 & -1 & -1 \end{array} \right ) $ the corresponding manifold is $X_{\Delta}=Bl_3(\mathbb{P}^2)$, the blow up of projective plane
at three $T$-equivariant points.

\bigskip

\hspace{-0.6cm} Consider the following three-fold examples:

\bigskip

\hspace{-0.6cm} \bf Example 2.3 \rm (Fano $\mathbb{P}^1$-bundles over $\mathbb{P}^2$): There are three Fano $\mathbb{P}^1$-bundles over $\mathbb{P}^2$ given by $X=\mathbb{P}(\mathcal{O}_{\mathbb{P}^2} \oplus \mathcal{O}_{\mathbb{P}^2}(k))$ for
$k=0,1,2$. The corresponding vertex matrix is $$ N_{\Delta^{\circ}}=\left ( \begin{array}{ccccc} 1 & 0 & 0 & -1 & 0 \\ 0 & 1 & 0 & -1 & 0 \\ 0 & 0 & 1 & k & -1 \end{array} \right ) $$

\bigskip

\hspace{-0.6cm} \bf Example 2.4 \rm (Fano $\mathbb{P}^2$-bundles over $\mathbb{P}^1$): There are two Fano $\mathbb{P}^2$-bundles over $\mathbb{P}^1$ given by $ X=\mathbb{ P}(\mathcal{ O} \oplus \mathcal{ O} \oplus \mathcal{ O}(k))$  for
$k=0,1$. The corresponding vertex matrix is $$ N_{\Delta^{\circ}}=\left ( \begin{array}{ccccc} 1 & 0 & 0 & -1 & -k \\ 0 & 1 & 0 & -1 & -k \\ 0 & 0 & 1 & 0 & -1 \end{array} \right ) $$

\bigskip

\hspace{-0.6cm} \bf Example 2.5 \rm (Fano $\mathbb{P}^1$-bundles over $\mathbb{P}^1 \times \mathbb{P}^1$): There are two Fano $\mathbb{P}^1$-bundles over $\mathbb{P}^1 \times \mathbb{P}^1$
given by $X=\mathbb{P}(\mathcal{ O}_{ \mathbb{ P}^1 \times \mathbb{ P}^1} \oplus  \mathcal{ O}_{\mathbb{ P}^1 \times \mathbb{ P}^1}(k_1,k_2))$ for $(k_1,k_2)=(0,0),(1,1),(1,-1)$. The vertex matrix is
$$ N_{\Delta^{\circ}} = \left ( \begin{array}{cccccc} 1 & 0 & 0 & k_1 & k_2 & -1 \\  0 & 1 & 0 & -1 & 0 & 0 \\ 0 & 0 & 1 & 0 & -1 & 0
\end{array} \right )$$

\bigskip

\hspace{-0.6cm} Denote by $\Delta(k)$ the set of $k$-dimensional faces of $\Delta$. As $T$ acts on $X$ it induces a decomposition of $X$ into orbits of the action. One of the fundamental properties of toric varieties is that $k$-dimensional orbits of the $T$-action are by themselves toric and are in one-to-one correspondence with
faces in $\Delta(k)$. Let $ V_X(F) \subset X$ be the closure of the orbit corresponding to the facet $F \in \Delta(k)$ in $X$. In particular, consider the group of toric divisors $$ Div_T(X) := \bigoplus_{F \in \Delta(n-1)} \mathbb{Z} \cdot V_X(F)$$ When $X$ is a smooth toric
manifold the group $Pic(X)$ admits a description in terms of the following short exact sequence $$ 0 \rightarrow M \rightarrow
Div_T(X) \rightarrow Pic(X) \rightarrow 0 $$ The map on the left hand side is given by $ m \rightarrow \sum_F \left < m, n_F \right > \cdot V_X(F) $ where $n_F \in \mathbb{N}_{\mathbb{R}}$ is the unit normal to the hyperplane spanned by the facet $F \in \Delta(n-1)$. In particular, note that $$ \rho(X)= rank \left ( Pic (X) \right )  = \vert \Delta(n-1) \vert -n $$ A toric divisor $D = \sum_F a_F \cdot V_X(F) \in Div_T(X) $ is associated with the following (possibly empty) polytope $$ \Delta_D:= \left \{ u \mid \left <u,n_F \right > \geq -a_F \textrm{ for all } F \right \} \subset M_{\mathbb{R}} $$ On the other, denote by $\mathcal{O}_X(D) \in Pic(X)$ the associated line bundle of $D$. The space of sections of $\mathcal{O}_X(D)$ is given in terms of the polytope $\Delta_D$ as follows $$H^0(X, \mathcal{O}_X(D)) \simeq \bigoplus_{m \in \Delta_D \cap M} m \cdot \mathbb{C}$$

\section{Derived Categories and Exceptional Collections}
\label{s:Dc}

\hspace{-0.6cm} Let $X$ be a smooth projective variety and let $\mathcal{D}^b(X)$ be the derived category of bounded complexes of
coherent sheaves of $\mathcal{O}_X$-modules. For a finite dimensional algebra $A$ denote by $\mathcal{D}^b(A)$ the derived category of bounded complexes of finite dimensional right modules over $A$. Given an object $T \in \mathcal{D}^b(X)$ denote by $A_T=Hom(T,T)$ the corresponding endomorphism algebra.

\bigskip

\hspace{-0.6cm}  \bf Definition 3.1: \rm An object $T \in \mathcal{D}^b(X)$ is called a \emph{tilting object} if the corresponding adjoint functors $$\begin{array}{ccc}  R Hom_X(T, -) : \mathcal{D}^b(X) \rightarrow \mathcal{D}^b(A_T) & ; & - \otimes^L_{A_T} T : \mathcal{D}^b(A_T) \rightarrow \mathcal{D}^b(X) \end{array} $$ are equivalences of categories. A locally free tilting object is called a tilting bundle.

\bigskip

\hspace{-0.6cm} One has the following characterization of tilting objects:

\bigskip

\hspace{-0.6cm} \bf Theorem 3.2 \rm (\cite{Bo,K}): An object $T \in \mathcal{D}^b(X)$ is a tilting object if and only if it satisfies
the following conditions:

\bigskip

- The endomorphism algebra $A_T$ is finite dimensional.

\bigskip

 - $Ext^i(T,T)=0$ for $0<i$.

\bigskip

 - The direct summands of $T$ generate $\mathcal{D}^b(X)$ as a triangulated category.

\bigskip

\hspace{-0.6cm} An object $ E \in \mathcal{D}^b(X)$ is said to be \emph{exceptional} if $Hom(E,E)=\mathbb{C}$ and $Ext^i(E,E)=0$ for $0<i$. An ordered collection $\left \{ E_1,...,E_N \right \} \subset \mathcal{D}^b(X)$
is said to be an \emph{exceptional collection} if each $E_j$ is exceptional and $$Ext^i(E_k,E_j) =0 \textrm{ for } j<k \textrm{ and } 0 \leq i $$ An
 exceptional collection is said to be \emph{strongly exceptional} if also $Ext^i(E_j,E_k)=0$ for $j \leq k$ and $0<i$. A strongly exceptional collection is called \emph{full} if its elements generate $\mathcal{D}^b(X)$ as a triangulated category. The importance of
full strongly exceptional collections in tilting theory is due to the following properties:

\bigskip

\hspace{-0.6cm} \bf Proposition 3.3 \rm (\cite{Bo,K}): Let $\mathcal{E}=\left \{ E_1,...,E_N \right \} \subset \mathcal{D}^b(X)$ be a collection.

\bigskip

(a) If $ \mathcal{E}$ is a full strongly exceptional collection then $T = \bigoplus_{i=1}^N E_i$ is a tilting object.

\bigskip

(b) If $T = \bigoplus_{i=1}^N E_i$ is a tilting object and $ \mathcal{E} \subset Pic(X) $ then $\mathcal{E}$ is a full strongly

\hspace{0.5cm} exceptional collection of line bundles.

\bigskip

\hspace{-0.6cm} Consider the following examples:

\bigskip

\hspace{-0.6cm} \bf Example 3.4 \rm (Projective space): \rm For $X = \mathbb{P}^n$ one has $Pic(X) \simeq H \cdot \mathbb{Z}$ where $H$ is the class represented by the normal bundle of a hyperplane in $X$.
The collection $$\mathcal{E}=\left \{ dH \mid d = 0,...,n \right \} \subset Pic(X)$$ is a full strongly exceptional collection, see \cite{B}.

\bigskip

\hspace{-0.6cm} \bf Example 3.5 \rm (Products): Let $X_1$ and $X_2$ be two projective manifolds and let $\mathcal{E}_1 \subset Pic(X_1)$ and $\mathcal{E}_2 \subset Pic(X_2)$ be full strongly exceptional collections on
$X_1$ and $X_2$ respectively. Then $$\mathcal{E}_1 \otimes \mathcal{E}_2 = \left \{ pr_1^{\ast}(L_1) \otimes pr_2^{\ast}(L_2) \vert L_1 \in \mathcal{E}_1 \textrm{ and } L_2 \in \mathcal{E}_2 \right \} \subset Pic(X_1 \times X_2)$$ is a full strongly exceptional collection on $X_1 \times X_2$.

\bigskip

\hspace{-0.6cm} \bf Example 3.6 \rm (Toric Del Pezzo surfaces): Exceptional collections for $\mathbb{P}^2$ and $\mathbb{P}^1 \times \mathbb{P}^1$ are given in examples 3.4 and 3.5. Recall that $$ Pic(Bl_k(\mathbb{P}^2)) \simeq H \cdot \mathbb{Z} \oplus \left ( \bigoplus_{i=1}^k E_i \cdot \mathbb{Z} \right ) $$ where $ E_i$ is the class of the normal bundle of the $i$-th exceptional divisor. It is shown in \cite{K} that $$ \mathcal{E} = \left \{ 0,H,H-E_i,2H-\sum_{i=1}^k E_i \vert i=1,...,k \right \} \subset Pic(Bl_k(\mathbb{P}^2))$$ is a full strongly exceptional collection on $Bl_k(\mathbb{P}^2)$ for $k=1,2,3$.

\bigskip

\hspace{-0.6cm} \bf Example 3.7 \rm (Threefold projective bundles): Let $ \pi : V \rightarrow B$ be a holomorphic vector bundle of rank $n-dim_{\mathbb{C}}(B)+1$ and let $X = \mathbb{P}(V)$. Express $Pic(X) \simeq \pi^{\ast} Pic(B) \oplus \xi \cdot \mathbb{Z}$ where $\xi$ is the tautological line bundle of $X$. Full strongly exceptional
collections on projective bundles were studied by Costa and Mir$\acute{\textrm{o}}$-Roig in \cite{CMR3} where they prove the following "key lemma": Let $\mathcal{E}_B = \left \{E_0,...,E_N \right \} \subset Pic(B)$ be a full strongly exceptional collection on $B$. Denote by $S^aV$ the $a$-th symmetric power
of $V$ and assume that $Hom(S^a \otimes E_m,E_l)=0$ for any $0 \leq a \leq n-dim_{\mathbb{C}}(B)$ and $0 \leq l \leq m \leq N$. Then $$\mathcal{E} = \left \{
\pi^{\ast} E_i \otimes k \xi \mid 0 \leq i \leq N, 0 \leq k \leq rank(V) \right \} \subset Pic(X)$$ is a full strongly exceptional collection on $X$. In particular, for examples 2.3-5, we have
$$ \mathcal{E} = \left \{ 0,H,2H,\xi,H+\xi,2H+\xi \right \} \subset Pic(X)$$ is a full strongly exceptional collection for $X=\mathbb{P}(\mathcal{O}_{\mathbb{P}^2} \oplus \mathcal{O}_{\mathbb{P}^2}(k))$ with $k=0,1,2$,
$$ \mathcal{E} = \left \{ 0,H,F,H+\xi,2 \xi,H+2\xi \right \} \subset Pic(X)$$ is a full strongly exceptional collection for $ X=\mathbb{ P}(\mathcal{ O} \oplus \mathcal{ O} \oplus \mathcal{ O}(k))$  with
$k=0,1$, and $$ \mathcal{E} = \left \{0,H_1,H_2,H_1+H_2,\xi,H_1+\xi,H_2+\xi,H_1+H_2+\xi \right \} \subset Pic(X)$$ is a full strongly exceptional collection for $X=\mathbb{P}(\mathcal{ O}_{ \mathbb{ P}^1 \times \mathbb{ P}^1} \oplus  \mathcal{ O}_{\mathbb{ P}^1 \times \mathbb{ P}^1}(k_1,k_2))$ with $(k_1,k_2)=(0,0),(1,1),(1,-1)$, see \cite{CMR3}.

\bigskip

\hspace{-0.6cm} \bf Remark 3.8: \rm Note that $X=Bl_1(\mathbb{P}^2)$ can be expressed as $X=\mathbb{P}(\mathcal{O}_{\mathbb{P}^1} \oplus \mathcal{O}_{\mathbb{P}^1}(1))$. In particular, one has $ \mathcal{E} = \left \{
0,H,H-E,2H-E \right\} = \left \{ 0,\xi,\pi^{\ast}H_{\mathbb{P}^1},\pi^{\ast}H_{\mathbb{P}^1}+\xi \right \} \subset Pic(X)$

\bigskip

\hspace{-0.6cm} \bf Remark 3.9 \rm (The Frobenius splitting method): A toric manifold $X$ is associated with a collection of maps $F_m : X \rightarrow X$ for $m \in \mathbb{N}$, known as Frobenius morphisms. It was shown by Thomsen \cite{T} that the push-forward
$(F_m)_{\ast}(\mathcal{O}_X)$ is a vector bundle, which splits as a sum of line bundles $\mathcal{D}_m \subset Pic(X)$, see also \cite{A}. The set $\mathcal{D}_m$ is independent of $m$ for $m >>0$ big enough, and we refer to the resulting "limit" collection
 $\mathcal{D}_X \subset Pic(X)$ as
the Frobenius collection of $X$. It is conjectured by Bondal that $\mathcal{D}_X$ generates $\mathcal{D}^b(X)$ as a triangulated category, see \cite{Bo3}.

\hspace{-0.6cm} When $\vert \mathcal{D}_X \vert = \rho(X)$, the Frobenius collection thus becomes a candidate for a full strongly exceptional collection. For instance, the full strongly exceptional collections $\mathcal{E}$ described in examples 3.6-8 are, in fact, Frobenius collections. Many further examples of toric manifolds whose Frobenius collections are full strongly exceptional collections were found, see \cite{CMR2,U}.

\hspace{-0.6cm} On the other hand, the condition $\vert \mathcal{D}_X \vert = \rho(X)$ does not always hold. For instance, when $X$ is a toric Fano threefold $\vert \mathcal{D}_X \vert = \rho(X)$ for sixteen of the eighteen Fano toric threefolds. Let us note that it is shown by Uehara in \cite{U} that in the cases when $\vert \mathcal{D}_X \vert =
\rho(X)$, the Frobenius collection $\mathcal{D}_X \subset Pic(X)$ is a full strongly exceptional collection of line bundles. Moreover, Uehara shows that, in
the remaining two cases, there is a subset $\mathcal{E} \subset \mathcal{D}_X$ which is a
full strongly exceptional collection.

\bigskip

\hspace{-0.6cm} A fundamental feature of tilting theory is the relation to quiver representations. Recall that a \emph{quiver} $Q=(Q_0,Q_1)$ is a directed graph such that $Q_0$ is the set of vertices of $Q$ and $Q_1$ is the set of directed edges of $Q$. Denote by $s,t : Q_1 \rightarrow Q_0$ the maps specifying the starting point $s(a) \in Q_0$ and end point $t(a) \in Q_0$ of an edge $a \in Q_1$. Denote by $Q_1(z,z') \subset Q_1$ the set of edges $a \in Q_1$ such that $s(a)=z$ and $t(a)=z'$.

\hspace{-0.6cm} A path in a quiver $Q$ is a sequence $a=a_1 \cdot ... \cdot a_n$ of edges such that $t(a_i)=s(a_{i+1})$. In particular, let $\mathbb{C} Q$ be the vector space spanned by all paths of $Q$ which admits the algebra operation given by $$ (a_1 \cdot ... \cdot a_n) \cdot (a_1' \cdot ... \cdot a_{n'}') =  \left  \{ \begin{array}{cc} a_1 \cdot ... \cdot a_n \cdot a_1' \cdot ... \cdot a_{n'}'
& t(a_n) = s(a'_1) \\ 0 & \textrm{ otherwise } \end{array} \right.  $$ we refer to the algebra $\mathbb{C}Q$ as the \emph{path algebra of the quiver} $Q$. A \emph{quiver with relations} is a pair $\widetilde{Q}= ( Q, R)$ where $Q$ is a quiver and $R \subset \mathbb{C}Q$ is a two sided ideal. In particular, we refer to $A_{\widetilde{Q}}=\mathbb{C}Q/R$ as the path algebra of the quiver with relations $\widetilde{Q}$. A finite dimensional right module over $A_{\widetilde{Q}}$ is called a representation of the quiver $\widetilde{Q}$.

\hspace{-0.6cm} Let $\mathcal{E} = \left \{ L_1,...,L_N \right  \} \subset Pic(X)$ be a collection and let $A_{\mathcal{E}} = \bigoplus_{i,j=1}^N Hom(L_i,L_j)$ be the corresponding endomorphism algebra. A choice of basis $$ B_{i,j} = \left \{ a_{i,j}^r \right \} \subset Hom(L_i,L_j)
\hspace{0.5cm} for  \hspace{0.3cm} 1 \leq i,j \leq N $$ determines a basis $B= \bigcup B_{i,j}$ of $A_{\mathcal{E}}$. A pair $(\mathcal{E},B)$ is associated with a quiver with relations $\widetilde{Q}(\mathcal{E},B)$ as follows: Let
 $Q_0(\mathcal{E},B)=\mathcal{E}$ be the set of vertices and let $Q_1(\mathcal{E},B;L_i,L_j)=B_{i,j}$ be the set of edges between $L_i,L_j \in \mathcal{E}$ for $1 \leq i,j \leq N$. In particular, a path $a=a_1 \cdot ... \cdot a_n$ in this (free) quiver can be considered as an element of $Hom(s(a_1),t(a_n))$. Thus, one has an exact sequence $$ 0 \rightarrow R(\mathcal{E},B) \rightarrow
\mathbb{C}Q(\mathcal{E},B) \rightarrow A_{\mathcal{E}} \rightarrow 0 $$ and taking the ideal of relations to be $R(\mathcal{E},B)$ gives $A_{\mathcal{E}}=A_{\widetilde{Q}(\mathcal{E},B)}$.

\hspace{-0.6cm} From now on we assume that $X$ is a toric manifold. Let $L \in Pic(X)$ be a line bundle, consider the set
$$ B(L) :=\left \{ D=\sum a_F \cdot V_X(F) \mid a_F \geq 0 \textrm{ and }
\mathcal{O}_X(D)=L \right \} \subset Div_T(X)$$ and denote by $i : B(L) \rightarrow H^0(X,L)$ the corresponding injection map. We refer to a line bundle $L \in Pic(X)$ as \emph{special} if
$i(B(L))$ is a basis for $H^0(X,L)$. We refer to a collection of line bundles $\mathcal{E} = \left \{L_1,...,L_N \right \} \subset Pic(X)$ as \emph{special} if $L_i-L_j$ is special for any $1 \leq i,j \leq N$. We observe the following:

\bigskip

\hspace{-0.6cm} \bf Proposition 3.10: \rm  For $X$ as in Example 2.2-5 the full strongly exceptional collection $\mathcal{E} = \mathcal{D}_X \subset Pic(X)$ is special.

\bigskip

\hspace{-0.6cm} Assume $\mathcal{E}$ is a special collection and let $B_{\mathcal{E}}: = \bigcup B(L_i-L_j)$ be the corresponding basis of $A_{\mathcal{E}}$. We refer to $\widetilde{Q}(\mathcal{E})=\widetilde{Q}(\mathcal{E},B_{\mathcal{E}})$ as the associated quiver of the special collection $\mathcal{E}$. Note that there exists a natural map $ D : \widetilde{Q}_1(\mathcal{E}) \rightarrow Div_T(X)$. We refer to the image $Div(\mathcal{E}) = D(\widetilde{Q}_1(\mathcal{E})) \subset Div_T(X)$ as the divisor set of the collection $\mathcal{E}$. For any $D \in Div(\mathcal{E})$ we denote by $\widetilde{Q}^D(\mathcal{E}) \subset \widetilde{Q}(\mathcal{E})$
the sub-quiver whose edge set is given by $\widetilde{Q}^D_1(\mathcal{E}):=\left \{ a \vert D(a) =D \right \} \subset \widetilde{Q}_1(\mathcal{E})  $. The following examples
describe the quiver $\widetilde{Q}(\mathcal{E})$ corresponding to the full strongly exceptional collections of Example 3.6-8:
\bigskip

\hspace{-0.6cm} \bf Example 3.11 \rm ($\widetilde{Q}(\mathcal{E})$ for toric Del-Pezzo manifolds):

\bigskip

\hspace{-0.6cm} (1) For $ X=\mathbb{ P}^2$ the vertex set $\Delta^{\circ}(0)$ is given by $$
\begin{array}{ccccc} n_1=(1,0) & ; & n_2=(0,1) & ; & n_3=(-1,1) \end{array} $$ with $[V_X(n_1)]=[V_X(n_2)]=[V_X(n_3)]=H$.
The quiver $\widetilde{Q}(\mathcal{E})$ for the exceptional collection $\mathcal{E}=\left \{0,H,2H \right \} \subset Pic(X) $ is thus given by $$

$$ with $$ \begin{array}{ccccc} D(a_1)=D(a_2) = V_X(n_1) & ; & D(b_1)=D(b_2) = V_X(n_2) & ; & D(c_1)=D(c_2) = V_X(n_3) \end{array}$$

\bigskip

\hspace{-0.6cm} (2) For $ X=\mathbb{ P}^1 \times \mathbb{ P}^1 $ the vertex set $\Delta^{\circ}(0)$ is given by $$ \begin{array}{ccccccc} n_1=(1,0) & ; & n_2=(0,1) & ; & n_3=(-1,0) & ; & n_4=(0,-1) \end{array} $$ with $[V_X(n_1)]=[V_X(n_3)]=H_1$ and  $[V_X(n_2)]=[V_X(n_4)]=H_2$.
The quiver $\widetilde{Q}(\mathcal{E})$ for the exceptional collection $\mathcal{E}=\left \{0,H_1,H_2,H_1+H_2 \right \} \subset Pic(X) $ is thus given by $$

$$ with $$ \begin{array}{cccc} D(a_1)=D(a_2) = V_X(n_1) & ; & D(b_1)=D(b_2) = V_X(n_2) & ;  \end{array}$$
$$ \begin{array}{cccc} D(c_1)=D(c_2) = V_X(n_3) & ; & D(d_1)=D(d_2)=V_X(n_4) & ; \end{array}$$
\bigskip

\hspace{-0.6cm} (3) For $X=Bl_1(\mathbb{ P}^2):$ the vertex set $\Delta^{\circ}(0)$ is given by $$ \begin{array}{ccccccc} n_1=(1,0) & ; & n_2=(0,1) & ; & n_3=(-1,0) & ; & n_4=(-1,-1) \end{array} $$ with $$ \begin{array}{ccccc}
[V_X(n_1)]=H=\xi & ; & [V_X(n_2)]=[V_X(n_4)]=H-E =\pi^{\ast}H_{\mathbb{P}^1} & ; & [V_X(n_3)]=E = \xi - \pi^{\ast}H_{\mathbb{P}^1} \end{array}$$ The quiver $\widetilde{Q}(\mathcal{E})$ for the exceptional collection $$ \mathcal{E}=\left \{0,H-E,H,2H-E \right \} = \left \{ 0, \pi^{\ast}H_{\mathbb{P}^1},\xi ,\xi +\pi^{\ast}H_{\mathbb{P}^1} \right \} \subset Pic(X) $$ is given by $$

$$ with $$ \begin{array}{cccccc} D(d_1)=D(d_2) = V_X(n_2) & ; & D(e_1)=D(e_2) = V_X(n_4) & ; &  D(f)=V_X(n_3) & ; \end{array}$$
$$ \begin{array}{cccc} D(a_1)=D(a_2) = V_X(n_2)+V_X(n_3) & ; & D(b_1)=D(b_2)=V_X(n_3)+V_X(n_4) & ; \end{array}$$
$$ \begin{array}{cc} D(c_1)=D(c_2) = V_X(n_1) & ; \end{array}$$

\bigskip

\hspace{-0.6cm} (4) For $ X=Bl_2(\mathbb{P}^2)$ the vertex set $\Delta^{\circ}(0)$ is given by $$ \begin{array}{ccccccccc} n_1=(1,0) & ; & n_2=(0,1) & ; & n_3=(-1,0) & ; & n_4=(-1,-1) & ; & n_5=(0,-1) \end{array} $$ with $$ \begin{array}{cccc}
[V_X(n_1)]=H-E_2 & ; & [V_X(n_2)]=H-E_1 & ;  \end{array}$$ $$ \begin{array}{ccccccc} & [V_X(n_3)]=E_1 & ; & [V_X(n_4)]=H-E_1-E_2 & ; & [V_X(n_5)]=E_2 & ; \end{array}$$ The quiver $\widetilde{Q}(\mathcal{E})$ for the exceptional collection $$ \mathcal{E}=\left \{0,H-E_1,H-E_2,H,2H-E_1-E_2 \right \} \subset Pic(X) $$ is thus given by  $$
$$

$$ \begin{array}{cccccc} D(e)= V_X(n_2)+V_X(n_3) & ; & D(f)= V_X(n_2)+V_X(n_4)+V_X(n_5) & ; &  D(g)=V_X(n_1)+V_X(n_5) & ; \end{array}$$
$$ \begin{array}{ccccc} D(a_1)= D(a_2) = V_X(n_2) & ; & D(b_1)=D(b_2)= V_X(n_4)+V_X(n_5) & ; \end{array}$$
$$ \begin{array}{ccccc} D(c_1)=D(c_2)= V_X(n_1) & ; & D(d_1)=D(d_2)= V_X(n_3)+V_X(n_4) & ;  \end{array}$$
$$ \begin{array}{cccccc} D(h)= V_X(n_3) & ; & D(i)=V_X(n_5) & ; & D(j)=V_X(n_4) \end{array}$$

\bigskip

\hspace{-0.6cm} (5) For $X=Bl_3(\mathbb{ P}^2)$ the vertex set $\Delta^{\circ}(0)$ is given by $$ \begin{array}{cccccc} n_1=(1,0) & ; & n_2=(1,1) & ; & n_3=(0,1) & ;\end{array} $$
$$ \begin{array}{cccccc} n_4=(-1,0) & ; & n_5=(-1,-1) & ; & n_6=(0,-1) & ; \end{array}$$ with $$ \begin{array}{cccccc}
[V_X(n_1)]=H-E_2-E_3 & ; & [V_X(n_2)]=E_3 & ; & [V_X(n_3)]=H-E_1-E_3 & ; \end{array}$$
$$ \begin{array}{cccccc} [V_X(n_4)]=E_1 & ; & [V_X(n_5)]=H-E_1-E_2 & ; & [V_X(n_6)]=E_2 & ; \end{array}$$
The quiver $\widetilde{Q}(\mathcal{E})$ for the exceptional collection $$ \mathcal{E}=\left \{0,H-E_1,H-E_2,H-E_3,H,2H-E_1-E_2-E_3 \right \} \subset Pic(X) $$ is given by  $$

$$ with $$ \begin{array}{cccc} D(a)= V_X(n_2)+V_X(n_3) & ; & D(b)= V_X(n_5)+V_X(n_6) & ; \end{array}$$
$$ \begin{array}{cccc} D(c)= V_X(n_1)+V_X(n_6) & ; & D(d)= V_X(n_3)+V_X(n_4) & ; \end{array}$$
$$ \begin{array}{cccc} D(e)= V_X(n_1)+V_X(n_2) & ; & D(f)= V_X(n_4)+V_X(n_5) & ; \end{array}$$
$$ \begin{array}{cccccc} D(g)= V_X(n_4) & ; & D(h)=V_X(n_2) & ; & D(i)=V_X(n_6) \end{array}$$
$$ \begin{array}{cccccc} D(j)= V_X(n_1) & ; & D(k)=V_X(n_5) & ; & D(l)=V_X(n_3) \end{array}$$

\bigskip

\hspace{-0.6cm} \bf Example 3.12 \rm ($\widetilde{Q}(\mathcal{E})$ for Fano $\mathbb{P}^1$-bundles over $\mathbb{P}^2$): For  $X=\mathbb{P}(\mathcal{O}_{\mathbb{P}^2} \oplus \mathcal{O}_{\mathbb{P}^2}(k))$ with $k=0,1,2$
 the vertex set $\Delta^{\circ}(0)$ is given by $$ \begin{array}{cccccc} n_1=(1,0,0) & ; & n_2=(0,1,0) & ; & n_3=(0,0,1) & ;\end{array} $$
$$ \begin{array}{cccc} n_4=(-1,-1,k) & ; & n_5=(0,0,-1) & ;  \end{array}$$ with $$ \begin{array}{cccccc}
[V_X(n_1)]= [V_X(n_2)]=[V_X(n_4) ]=\pi^{\ast} H & ; & [V_X(n_3)]=\xi-k\pi^{\ast}H & ; & [V_X(n_5)]=\xi & ; \end{array}$$
Consider the exceptional collection $$ \mathcal{E}=\left \{0,\xi,\pi{^\ast}H,\pi^{\ast}H+\xi,2\pi^{\ast}H,2\pi^{\ast}H+\xi \right \} \subset Pic(X) $$ (i) For $k=0$ the quiver $\widetilde{Q}(\mathcal{E})$ is given by

$$

$$ with
$$ \begin{array}{cccccc} D(a_i)=V_X(n_1) & ; & D(b_i)= V_X(n_2) & ; & D(c_i)=V_X(n_4) \end{array}$$
$$ \begin{array}{cccc} D(d_j)= V_X(n_3) & ; & D(e_j)= V_X(n_5) & ; \end{array}$$ where $ 1 \leq i \leq 4$ and $1 \leq j \leq 3$.

\bigskip

\hspace{-0.6cm} (ii) For $k=1$ the quiver is: $$
$$ with $$ \begin{array}{cccccc} D(a_i)=V_X(n_1) & ; & D(b_i)= V_X(n_2) & ; & D(c_i)=V_X(n_4) \end{array}$$
$$ \begin{array}{cccc} D(d_k)= V_X(n_3) & ; & D(e_j)= V_X(n_5) & ; \end{array}$$
 where $ 1 \leq i \leq 4$ and $1 \leq j \leq 3$ and $1 \leq k \leq 2$.

\bigskip

\hspace{-0.6cm} (iii) For $k=2$ the quiver is: $$$$
 with $$ \begin{array}{cccccc} D(a_i)=V_X(n_1) & ; & D(b_i)= V_X(n_2) & ; & D(c_i)=V_X(n_4) \end{array}$$
$$ \begin{array}{cccc} D(d)= V_X(n_3) & ; & D(e_j)= V_X(n_5) & ; \end{array}$$
 where $ 1 \leq i \leq 4$ and $1 \leq j \leq 3$.
\bigskip

 \hspace{-0.6cm} \bf Example 3.13 \rm ($\widetilde{Q}(\mathcal{E})$ for Fano $\mathbb{P}^2$-bundles over $\mathbb{P}^1$): For $ X=\mathbb{ P}(\mathcal{ O}_{\mathbb{P}^1} \oplus \mathcal{ O}_{\mathbb{P}^1} \oplus \mathcal{ O}_{\mathbb{P}^1} (k))$ with $k=0,1$
 the vertex set $\Delta^{\circ}(0)$ is given by $$ \begin{array}{cccccc} n_1=(1,0,0) & ; & n_2=(0,1,0) & ; & n_3=(0,0,1) & ;\end{array} $$
$$ \begin{array}{cccc} n_4=(-1,-1,0) & ; & n_5=(-k,-k,-1) & ;  \end{array}$$ with $$ \begin{array}{cccccc}
[V_X(n_1)]= [V_X(n_2)]= \xi & ; & [V_X(n_3)]=\xi-k\pi^{\ast}H & ; & [V_X(n_3)]=[V_X(n_5)]=\pi^{\ast}H & ; \end{array}$$
The quiver $\widetilde{Q}_{\mathcal{E}}$ for the exceptional collection $$ \mathcal{E}=\left \{0,\xi,2\xi,\pi^{\ast}H,\pi^{\ast}H+\xi,\pi^{\ast}H+2\xi \right \} \subset Pic(X) $$ is given by
$$

$$  with $$ \begin{array}{cccccc} D(a_i)=V_X(n_1) & ; & D(b_i)= V_X(n_2) & ; & D(c_j)=V_X(n_3) & ; \end{array} $$
$$ \begin{array}{cccc} D(d_j) = V_X(n_5) & ; & D(e_k) = V_X(n_4) & ; \end{array}$$ where $ 1 \leq i \leq 4$ and $1 \leq j \leq 3$ and $ 1 \leq k \leq 2$.

\bigskip

\hspace{-0.6cm} \bf Example 3.14 \rm ($\widetilde{Q}(\mathcal{E})$ for $\mathbb{P}^1$-bundles over $\mathbb{P}^1 \times \mathbb{P}^1$): For
$X=\mathbb{P}(\mathcal{ O}_{ \mathbb{ P}^1 \times \mathbb{ P}^1} \oplus  \mathcal{ O}_{\mathbb{ P}^1 \times \mathbb{ P}^1}(k_1,k_2))$ with $(k_1,k_2)=(0,0),(1,1),(1,-1)$
 the vertex set $\Delta^{\circ}(0)$ is given by $$ \begin{array}{cccccc} n_1=(1,0,0) & ; & n_2=(0,1,0) & ; & n_3=(0,0,1) & ;\end{array} $$
$$ \begin{array}{cccccc} n_4=(k_1,-1,0) & ; & n_5=(k_2,0,-1) & ; & n_6=(-1,0,0) & ;  \end{array}$$ with $$ \begin{array}{cccccc}
[V_X(n_1)]= \xi & ; & [V_X(n_2)]=[V_X(n_4)]=\pi^{\ast}H_1 & ; & [V_X(n_3)]=[V_X(n_5)]=\pi^{\ast}H_2 & ; \end{array}$$
$$ \begin{array}{cc} [V_X(n_6)]=\xi+k_1\pi^{\ast}H_1+k_2\pi^{\ast}H_2 & ; \end{array}$$
Consider the exceptional collection $$ \mathcal{E}=\left \{0,\xi,\pi^{\ast}H_1,\pi^{\ast}H_1+\xi,\pi^{\ast}H_2,\pi^{\ast}H_2+\xi,
\pi^{\ast}H_1+\pi^{\ast}H_2,\pi^{\ast}H_1+\pi^{\ast}H_2+\xi \right \} \subset Pic(X) $$ For $(k_1,k_2)=(0,0)$ the quiver $\widetilde{Q}(\mathcal{E})$ is given by

$$

$$  $$ \begin{array}{cccccc} D(a_i)=V_X(n_2) & ; & D(b_i)= V_X(n_4) & ; & D(c_i)=V_X(n_3) & ; \end{array} $$
 $$ \begin{array}{cccccc} D(d_i)=V_X(n_5) & ; & D(e_i)= V_X(n_1) & ; & D(f_i)=V_X(n_6) & ; \end{array} $$ where $1 \leq i \leq 4$.

\bigskip

\hspace{-0.6cm} (ii) For $(k_1,k_2)=(1,1)$ the quiver $\widetilde{Q}(\mathcal{E})$ is given by

$$

$$
 $$ \begin{array}{cccccc} D(a_i)=V_X(n_2) & ; & D(b_i)= V_X(n_4) & ; & D(c_i)=V_X(n_3) & ; \end{array} $$
 $$ \begin{array}{cccccc} D(d_i)=V_X(n_5) & ; & D(e_i)= V_X(n_1) & ; & D(f)=V_X(n_6) & ; \end{array} $$ where $1 \leq i \leq 4$.

\bigskip

\hspace{-0.6cm} (iii) For $(k_1,k_2)=(1,-1)$ the quiver $\widetilde{Q}(\mathcal{E})$ is given by
$$

$$
 $$ \begin{array}{cccccc} D(a_i)=V_X(n_2) & ; & D(b_i)= V_X(n_4) & ; & D(c_i)=V_X(n_3) & ; \end{array} $$
 $$ \begin{array}{cccccc} D(d_i)=V_X(n_5) & ; & D(e_i)= V_X(n_1) & ; & D(f)=V_X(n_6) & ; \end{array} $$ where $1 \leq i \leq 4$.

\section{The Landau-Ginzburg System and Monodromy}
\label{s:LGM}

\hspace{-0.6cm} Let $X$ be a $n$-dimensional toric Fano manifold given by a Fano polytope $\Delta \subset M_{\mathbb{R}}$ and let $\Delta^{\circ} \subset N_{\mathbb{R}}$ be the corresponding polar polytope. Let $L(\Delta^{\circ}) \subset \mathbb{C}[N]$ be the space of Laurent polynomials whose Newton polytope is $\Delta^{\circ}$ and let $ f_X = \sum_{n \in \Delta^{\circ}(0)} z^n \in L(\Delta^{\circ})$ be the Landau-Ginzburg potential of $X$. We refer to $$ z_i \frac{\partial}{\partial z_i } f_X(z_1,...,z_n)=0 \hspace{0.5cm} \textrm{ for } \hspace{0.25cm} i=1,...,n $$ as the Landau-Ginzburg system of equations of $X$ and denote by $Crit(f_X) \subset (\mathbb{C}^{\ast})^n$ the corresponding solution scheme. The Landau-Ginzburg potential was first introduced by Batyrev in \cite{Ba}, in the context of the study of the small quantum
cohomology $QH(X)$ of $X$. The main property of the Landau-Ginzburg potential is the existence of a ring isomorphism $QH(X) \simeq Jac(f_X)$ where $$ Jac(f_X) : = \frac{\mathbb{C}[N]} {\left < z_i \frac{\partial}{\partial z_i} f_X \mid i=1,...,n \right >}$$ is the \emph{Jacobian ring} of $f_X$, see \cite{Ba,OT,FOOO}. Ostrover and Tyomkin describe the following semi-simplicity criteria: $QH(X)$ is semi-simple if and only if $Crit(f_X)$ is a reduced scheme, see \cite{OT}. We refer to a toric Fano manifold as \emph{semi-simple} if its solution scheme $Crit(f_X)$ is reduced.

\hspace{-0.6cm} We observe that the solution scheme $Crit(f_X) \subset (\mathbb{C}^{\ast})^n$ of a semi-simple toric Fano manifold $X$ admits a natural mondromy action. Indeed, consider an element $ f \in L(\Delta^{\circ}) $ as a linear functional on $L(\Delta^{\circ})^{\vee}=Hom(L(\Delta^{\circ}),\mathbb{C})$. For a non-zero element $ f \in L(\Delta^{\circ})$ denote by $H(f) \subset \mathbb{P}(L(\Delta^{\circ})^{\vee})$ the hyperplane
given by the kernel of $f$.  Note that, if $f$ is non-constant, the hyperplane section $\Sigma(f):=X^{\circ} \cap H_f$ is given as the closure of the hypersurface $\left \{z \vert f(z) =0 \right \} \subset (\mathbb{C}^{\ast})^n $ in $X^{\circ}$. The definition of the solution shceme could hence be generalized for any non-constant $f \in L(\Delta^{\circ})$ by setting $$Crit(f) := \Sigma \left ( z_1 \frac{ \partial}{\partial z_1} f \right ) \cap ... \cap  \Sigma \left ( z_n \frac{ \partial}{\partial z_n} f \right )  \subset X^{\circ} $$ In particular, we refer to the hypersurface
$$ R_X := \left \{ f \mid Crit(f) \textrm{ is non-reduced} \right \} \subset L(\Delta^{\circ})$$ as the \emph{resultant hypersurface} of $X$. A path of the form $\gamma : [0,1] \rightarrow L(\Delta^{\circ}) \setminus R_X$ gives rise to a map $C_{\gamma} : Crit(\gamma(0)) \rightarrow Crit(\gamma(1))$ via analytic continuation. By standard considerations, the map $C_{\gamma}$ is an invariant of the homotopy class of $\gamma$ in $\pi_1(L(\Delta^{\circ}) \setminus R_X; \gamma(0),\gamma(1))$. In particular, assuming  $X$ is a semi-simple toric Fano manifold, we obtain the following monodromy action $$ M: \pi_1(L(\Delta^{\circ}) \setminus R_X; f_X) \rightarrow Aut(Crit(f_X))$$ given by $[\gamma] \mapsto C_{[\gamma]}$. For an element $[\gamma] \in \pi_1( L(\Delta^{\circ}) \setminus R_X , f_X)$ denote by $\widetilde{Q}([\gamma])$ the quiver whose vertex set is $\widetilde{Q}_0([\gamma]) = Crit(f_X)$, and has an edge $a_{[\gamma]}(z) \in
\widetilde{Q}_1([\gamma])$ such that $s(a_{[\gamma]}(z))=z$ and $t(a_{[\gamma]}(z))= M([\gamma])(z)$, for any $z \in Crit(f_X)$.

\section{ The Exceptional Maps $\bf L : Crit(f_X) \rightarrow Pic(X)$ and Monodromy}
\label{s:EM}

\hspace{-0.6cm} A map $L : Crit(f_X) \rightarrow Pic(X)$ is said to be exceptional if $\mathcal{E}_L := L(Crit(f_X)) \subset Pic(X)$ is a full strongly exceptional collection. The map $L$ is said to be special if $\mathcal{E}_L$ is a special collection, in the sense of Section 3.  For each of the manifolds (a)-(d) we associate a map of the form $\Gamma : Div_T(X) \rightarrow \pi_1(L(X) \setminus R_X,f_X)$.

\hspace{-0.6cm} For a toric divisor $D = \sum a_F V_X(F) \in Div_T(X)$ and an element $f(z) = \sum b_F z^{n_F} \in L(\Delta^{\circ})$ let $\gamma_{(f,D)} : S^1 \rightarrow L(\Delta^{\circ})$ be the loop given by $\gamma_{(f,D)}(t,z) = \sum b_F e^{ia_Ft} z^{n_F}$.
 
\hspace{-0.6cm} For the toric Del-Pezzo manifolds of type (a) we define the map $\Gamma$ by $\Gamma(D) := [\gamma_{(f_X,D)}]$. For the projective bundles of type (b)-(d) express the Landau-Ginzburg potential as $f_X(z) = f^{base}_X(z) + f^{fiber}_X(z)$. For $0 < 
\epsilon$ denote by $$g^t(z) =(1-t+ \epsilon t) f_X^{base} + f_X^{fiber} \in L(\Delta^{\circ})$$ for $t \in [0,1]$. Set $$ \widetilde{\gamma}_D(t,z) = \left \{ \begin{array}{cc} g^{3t}(z) & t \in [0,\frac{1}{3})  \\ \gamma_{(g^1,D)}(3t-1,z) & t \in [\frac{1}{3},\frac{2}{3}) \\ g^{3-3t}(z) & t \in [\frac{2}{3},1] \end{array} \right. $$ and set $\Gamma(D) := [\widetilde{\gamma}_D^{\epsilon}]$ for $0 < \epsilon $  small enough. Note that $\widetilde{\gamma}^1_D = \gamma_{(f_X,D)}$. In practice, $\epsilon=1$ is taken for all examples aside from $\mathbb{P}(\mathcal{O} \oplus \mathcal{O}(2))$ of class (b) and $\mathbb{P}(\mathcal{O}_{\mathbb{P}^1 \times \mathbb{P}^1} \oplus \mathcal{O}_{\mathbb{P}^1 \times \mathbb{P}^1}(1,-1) )$ of class (d) for which we take $\epsilon = \frac{1}{2}$. 

\hspace{-0.6cm} For the manifolds (a)-(d) denote by $\widetilde{Q}(D) = \widetilde{Q}(\Gamma(D))$. We say that a special exceptional map $L : Crit(f_X) \rightarrow Pic(X)$ is \emph{$M$-aligned} ($M$ stands for monodromy) if it satisfies the following condition $$ \widetilde{Q}^D(\mathcal{E}) \subset \widetilde{Q}([\gamma_D]) \hspace{0.5cm} \textrm{for any} \hspace{0.25cm} D \in Div(\mathcal{E})$$ By definition, any toric Fano manifold $X$ which admits a full strongly exceptional collection $\mathcal{E} \subset Pic(X)$, also admits an exceptional map. Note, however, that the existence of a $M$-aligned exceptional map is a far less trivial property. Our main result is:

\bigskip

\hspace{-0.6cm} \bf Theorem 5.1: \rm Let $X$ be as in Example 2.2-5. Then $X$ admits a $M$-aligned exceptional map $L : Crit(f_X) \rightarrow Pic(X)$.

\bigskip

\hspace{-0.6cm} The rest of this section is devoted to the definition of the exceptional maps and the verification of the $M$-aligned property.

\subsection{Definition of the Exceptional Maps $\bf L : Crit(f_X) \rightarrow Pic(X)$:}
\label{s:EMD} This sub-section is devoted to the definition the exceptional maps $L : Crit(f_X) \rightarrow Pic(X)$ for the toric Fano manifolds of Example 2.2-5. In particular, we give a numerical description of the solution set
$Crit(f_X) \subset (\mathbb{C}^{\ast})^n$ of these manifolds.

\bigskip

\hspace{-0.6cm}  5.1.1 \bf The Del Pezzo surface case \rm

\bigskip

\hspace{-0.6cm} (i) For $\underline{ X=\mathbb{ P}^2}$ the LG-potential is given by $f_X(z_1,z_2) = z_1+z_2+\frac{1}{z_1z_2}$ and the Landau-Ginzburg system is $$ \begin{array}{ccc} z_1 \frac{\partial}{\partial z_1} f_X( z_1,z_2) =
z_1 - \frac{1}{z_1 z_2} =0 & ; &  z_2 \frac{\partial}{\partial z_2} f_X( z_1,z_2) =
z_2 - \frac{1}{z_1 z_2} =0 \end{array} $$ The solution set $Crit(f_X) \subset (\mathbb{C}^{\ast})^2$ is given by $$ \begin{array}{cccccc} z^0=(1,1) & ; & z^1=(e^{\frac{2 \pi i}{3}},e^{\frac{2 \pi i}{3}}) & ; & z^2=(e^{\frac{4 \pi i}{3}},e^{\frac{4 \pi i}{3}})
& ; \end{array}$$
The following is a graphical description:
$$\begin{array}{ccc} 
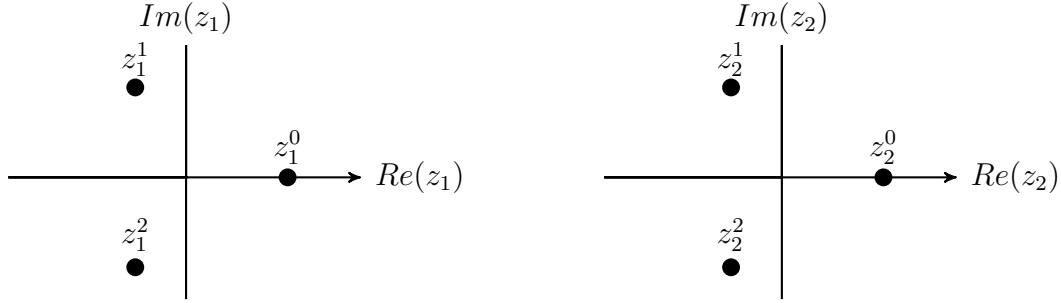
\begin{tikzpicture}[->,>=stealth',shorten >=1pt,auto,node distance=2.5cm,main node/.style={circle,draw,font=\sffamily\bfseries \small}, scale=1.35]
    % Draw axes
    \draw [->,thick] (0,-1.2)--(0,1.3) node (yaxis) [above] {$Im(z_1)$}
        |- (-1.75,0)--(1.75,0) node (xaxis) [right] {$Re(z_1)$};

    \node (a) at (1,0) {};
    \node (b) at (-0.5,0.886){};
    \node (c) at (-0.5,-0.886){} ;
  \node (a1) at (1,0.3) {$z^0_1$};
    \node (b1) at (-0.5,1.166) {$z^1_1$};
    \node (c1) at (-0.5,-0.566) {$z^2_1$};
   \fill[black] (a) circle (2.5pt);
\fill[black] (b) circle (2.5pt);
\fill[black] (c) circle (2.5pt);
\end{tikzpicture}

& \hspace{1cm}
&
\begin{tikzpicture}[->,>=stealth',shorten >=1pt,auto,node distance=2.5cm,main node/.style={circle,draw,font=\sffamily\bfseries \small}, scale=1.35]
    % Draw axes
    \draw [->,thick] (0,-1.2)--(0,1.3) node (yaxis) [above] {$Im(z_2)$}
        |- (-1.75,0)--(1.75,0) node (xaxis) [right] {$Re(z_2)$};

    \node (a) at (1,0){};
    \node (b) at (-0.5,0.886){};
    \node (c) at (-0.5,-0.886) {};
  \node (a1) at (1,0.3) {$z^0_2$};
    \node (b1) at (-0.5,1.166) {$z^1_2$};
    \node (c1) at (-0.5,-0.566) {$z^2_2$};
   \fill[black] (a) circle (2.5pt);
\fill[black] (b) circle (2.5pt);
\fill[black] (c) circle (2.5pt);

\end{tikzpicture}
\end{array}$$ we define the exceptional map $L : Crit(f_X) \rightarrow Pic(X)$ to be $L(z^k) = kH$ for $k=0,1,2$.

\bigskip

\hspace{-0.6cm} (ii) For $\underline{ X=\mathbb{ P}^1 \times \mathbb{ P}^1 }$ the LG-potential is given by $f_X(z_1,z_2) = z_1+z_2+\frac{1}{z_1}+\frac{1}{z_2}$ and the Landau-Ginzburg system is $$ \begin{array}{ccc} z_1 \frac{\partial}{\partial z_1} f_X( z_1,z_2) =
z_1 - \frac{1}{z_1 } =0 & ; &  z_2 \frac{\partial}{\partial z_2} f_X( z_1,z_2) =
z_2 - \frac{1}{z_2} =0 \end{array} $$ The solution set $Crit(f_X) \subset (\mathbb{C}^{\ast})^2$ is given by $$ \begin{array}{cccccccc} z^{00}=(1,1) & ; & z^{01}=(1,-1) & ; & z^{10}=(-1,1) & ; & z^{11}=(-1,1) & ; \end{array}$$
The following is a graphical description:

$$ 
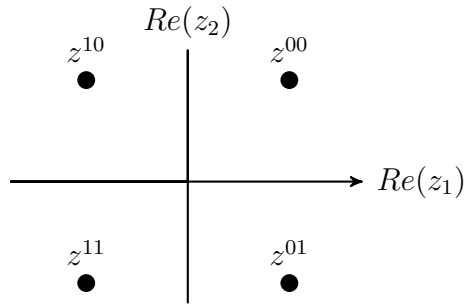
\begin{tikzpicture}[->,>=stealth',shorten >=1pt,auto,node distance=2.5cm,main node/.style={circle,draw,font=\sffamily\bfseries \small}, scale=1.35]
    % Draw axes
    \draw [->,thick] (0,-1.2)--(0,1.3) node (yaxis) [above] {$Re(z_2)$}
        |- (-1.75,0)--(1.75,0) node (xaxis) [right] {$Re(z_1)$};

    \node (a) at (1,1){};
    \node (b) at (-1,1){};
    \node (c) at (1,-1){} ;
\node (d) at (-1,-1){};
    \node (a1) at (1,1.3) {$z^{00}$};
    \node (b1) at (-1,1.3) {$z^{10}$};
    \node (c1) at (1,-0.7) {$z^{01}$};
\node (d1) at (-1,-0.7) {$z^{11}$};
   \fill[black] (a) circle (2.5pt);
\fill[black] (b) circle (2.5pt);
\fill[black] (c) circle (2.5pt);
\fill[black] (d) circle (2.5pt);
\end{tikzpicture}
$$ we define the exceptional map $L : Crit(f_X) \rightarrow Pic(X)$ to be $L(z^{k_1,k_2}) = k_1H_1+k_2H_2$ for $0 \leq k_1,k_2 \leq 1$.

\bigskip

\hspace{-0.6cm} (iii) $\underline{ X=Bl_1(\mathbb{ P}^2):}$ The LG-potential is given by $f_X(z_1,z_2) = z_1+z_2+\frac{1}{z_1z_2}+\frac{1}{z_1}$ and the Landau-Ginzburg system is $$ \begin{array}{ccc} z_1 \frac{\partial}{\partial z_1} f_X( z_1,z_2) =
z_1 - \frac{1}{z_1 z_2}-\frac{1}{z_1} =0 & ; &  z_2 \frac{\partial}{\partial z_2} f_X( z_1,z_2) =
z_2 - \frac{1}{ z_2} =0 \end{array} $$ The following is a numerical approximation of the elements of $Crit(f_X) \subset (\mathbb{C}^{\ast})^2$: $$ \begin{array}{cccc} z^0 \approx (1.49,0.81) & ; & w^1 \approx (0.52,-1.38) & ; \end{array}$$
$$ \begin{array}{cccc} z^1 \approx (-1+0.51i,-0.21+0.91i) & ; & z^2 \approx (-1-0.51i,-0.21-0.91i)& ; \end{array}$$
The following is a graphical description:
$$\begin{array}{ccc} 
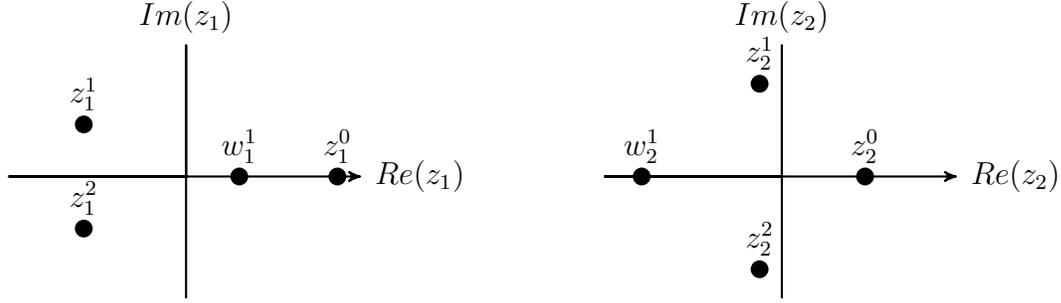
\begin{tikzpicture}[->,>=stealth',shorten >=1pt,auto,node distance=2.5cm,main node/.style={circle,draw,font=\sffamily\bfseries \small}, scale=1.35]
    % Draw axes
    \draw [->,thick] (0,-1.2)--(0,1.3) node (yaxis) [above] {$Im(z_1)$}
        |- (-1.75,0)--(1.75,0) node (xaxis) [right] {$Re(z_1)$};

    \node (a) at (0.524,0){} ;
    \node (b) at (-1.007,0.513){};
    \node (c) at (-1.007,-0.513){} ;
    \node (d) at (1.490,0){} ;
\node (a1) at (0.524,0.3) {$w_1^1$};
    \node (b1) at (-1.007,0.813) {$z^1_1$};
    \node (c1) at (-1.007,-0.213) {$z^2_1$} ;
    \node (d1) at (1.490,0.3) {$z_1^0$};
   \fill[black] (a) circle (2.5pt);
\fill[black] (b) circle (2.5pt);
\fill[black] (c) circle (2.5pt);
\fill[black] (d) circle (2.5pt);
\end{tikzpicture}

& \hspace{1cm}
&
\begin{tikzpicture}[->,>=stealth',shorten >=1pt,auto,node distance=2.5cm,main node/.style={circle,draw,font=\sffamily\bfseries \small}, scale=1.35]
    % Draw axes
    \draw [->,thick] (0,-1.2)--(0,1.3) node (yaxis) [above] {$Im(z_2)$}
        |- (-1.75,0)--(1.75,0) node (xaxis) [right] {$Re(z_2)$};

    \node (a) at (-1.380,0){};
    \node (b) at (-0.219,-0.914){};
    \node (c) at (-0.219,0.914) {};
\node (d) at (0.819,0) {};
   \node (a1) at (-1.380,0.3) {$w^1_2$};
    \node (b1) at (-0.219,-0.614) {$z^2_2$};
    \node (c1) at (-0.219,1.214) {$z^1_2$};
\node (d1) at (0.819,0.3) {$z^0_2$};
   \fill[black] (a) circle (2.5pt);
\fill[black] (b) circle (2.5pt);
\fill[black] (c) circle (2.5pt);
\fill[black] (d) circle (2.5pt);
\end{tikzpicture}
\end{array}$$ Note that contrary to the previous examples, in this case, the definition of the exceptional map is not directly evident from the numerical data. However, the exceptional map is "uncovered"
by utilizing the fact that $X$ is given geometrically as the blow up of $\mathbb{P}^2$ at one point.

\hspace{-0.6cm} Indeed, consider the following one-parametric family of elements $$g^t(z_1,z_2)=z_1+z_2+\frac{1}{z_1z_2} + \frac{1-t}{z_1} \in L(\Delta^{\circ}) $$
for $t = [0,1)$. Note that $g^0(z_1,z_2)=f_X(z_1,z_2)$ and $g^t(z_1,z_2)  \rightarrow f_{\mathbb{P}^2}(z_1,z_2)$ for $t \rightarrow 1^-$. By analytic continuation one can express $Crit(g^t) = \left \{ z^0(t),z^1(t),z^2(t),w^1(t) \right \} \subset (\mathbb{C}^{\ast})^2$.
 Direct computation shows $$ \begin{array}{cccc} z^0(t) \rightarrow z^0_{\mathbb{P}^2}=(1,1) & ; & w^1(t) \rightarrow w^1(1)=(\infty,0) & ; \end{array}$$
$$ \begin{array}{cccc} z^1(t) \rightarrow z^1_{\mathbb{P}^2}=(e^{\frac{2 \pi i}{3}},e^{\frac{2 \pi i}{3}})  & ; & z^2(t) \rightarrow z^2_{\mathbb{P}^2}=(e^{\frac{4 \pi i}{3}},e^{\frac{4 \pi i}{3}}) & ; \end{array}$$ as $t \rightarrow 1^-$. In view of this we define
the exceptional map $L : Crit(f_X) \rightarrow Pic(X)$ by $$ \begin{array}{cccccccc} L(z^0)=0 & ; & L(z^1)=H & ; & L(z^2)=2H-E & ; & L(w^1)=H-E & ; \end{array}$$

\hspace{-0.6cm} Note that, viewing $X$ as a $\mathbb{P}^1$-bundle over $\mathbb{P}^1$, one can define the one-parametric family of elements $$h^t(z_1,z_2)=z_1+(1-t)z_2+\frac{1-t}{z_1z_2} + \frac{1}{z_1} \in L(\Delta^{\circ})$$
for $t = [0,1)$. Note that $h^0(z_1,z_2)=f_X(z_1,z_2)$ and $h^t(z_1,z_2)  \rightarrow f_{\mathbb{P}^1}(z_1,z_2)$ for $t \rightarrow 1^-$. Express $Crit(h^t) = \left \{ \widetilde{z}^0(t),\widetilde{z}^1(t),\widetilde{z}^2(t),\widetilde{w}^1(t) \right \} \subset (\mathbb{C}^{\ast})^2$. Direct computation shows
$$ \begin{array}{cccc} \widetilde{z}^0(t),\widetilde{w}^1(t) \rightarrow \widetilde{z}^0_{\mathbb{P}^1}=(1,0)  & ; & \widetilde{z}^1(t),\widetilde{z}^2(t) \rightarrow \widetilde{z}^1_{\mathbb{P}^1}=(-1,0) & ; \end{array}$$ when $t \rightarrow 1$. Which leads to the
analogous definition of $L : Crit(f_X) \rightarrow Pic(X)$ as $$ \begin{array}{cccccccc} L(z^0)=0 & ; & L(z^1)=\xi & ; & L(z^2)=\xi+\pi^{\ast}H & ; & L(w^1)=\pi^{\ast}H & ; \end{array}$$

%$$ \begin{tikzpicture}[scale=1.25,->,>=stealth',shorten >=1pt,auto,node distance=2.5cm,main node/.style={circle,draw,font=\sffamily\bfseries },main node1/.style={circle,fill=black,draw,font=\sffamily\bfseries }]
%
%    \node[main node] (a) at (0.524,0);
%    \node[main node] (b) at (-1.007,0.513);
%    \node[main node] (c) at (-1.007,-0.513) ;
%    \node[main node] (d) at (1.490,0) ;
% \node[main node1] (e) at (1,0);
 %   \node[main node1] (f) at (-0.5,0.866);
%    \node[main node1] (g) at (-0.5,-0.866) ;
%  \node[main node1] (h) at (0,0) ;
% \path[every node/.style={font=\sffamily\small}]
%(a) edge (h)
%(b) edge (f)
%(c) edge (g)
%(d) edge (e);
%
%\end{tikzpicture} $$

\bigskip

\hspace{-0.6cm} (iv) For $\underline{X=Bl_2(\mathbb{ P}^2)}$ the LG-potential is given by $f_X(z_1,z_2) = z_1+z_2+\frac{1}{z_1z_2}+\frac{1}{z_1}+\frac{1}{z_2}$ and the Landau-Ginzburg system is $$ \begin{array}{ccc} z_1 \frac{\partial}{\partial z_1} f_X( z_1,z_2) =
z_1 - \frac{1}{z_1 z_2} -\frac{1}{z_1}=0 & ; &  z_2 \frac{\partial}{\partial z_2} f_X( z_1,z_2) =
z_2 - \frac{1}{z_1 z_2} -\frac{1}{z_2}=0 \end{array} $$ The following is a numerical approximation of the elements of $Crit(f_X) \subset (\mathbb{C}^{\ast})^2$: $$ \begin{array}{cccccc} z^0 \approx (1.32,1.32) & ; &
w^1 \approx (0.61,-1.61) & ; & w^2 \approx (-1.61,0.61) & ;  \end{array}$$
$$ \begin{array}{cccc} z^1 \approx (-0.66+0.56i,-0.66+0.56i) & ; & z^2 \approx (-0.66-0.56i,-0.66-0.56i)& ; \end{array}$$
The following is a graphical description:
$$\begin{array}{ccc} 
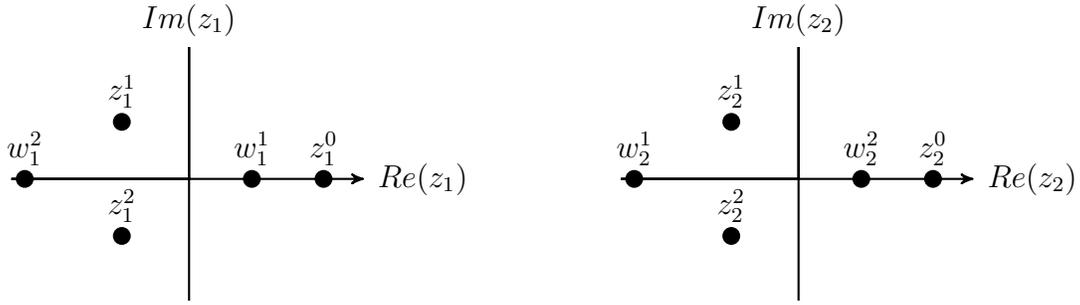
\begin{tikzpicture}[->,>=stealth',shorten >=1pt,auto,node distance=2.5cm,main node/.style={circle,draw,font=\sffamily\bfseries \small}, scale=1.35]
    % Draw axes
    \draw [->,thick] (0,-1.2)--(0,1.3) node (yaxis) [above] {$Im(z_1)$}
        |- (-1.75,0)--(1.75,0) node (xaxis) [right] {$Re(z_1)$};

    \node (a) at (-1.618,0){};
    \node (b) at (0.618,0){};
    \node (c) at (1.324,0) {};
    \node (d) at (-0.662,0.562) {};
  \node (e) at (-0.662,-0.562) {};
  \node (a1) at (-1.618,0.3) {$w^2_1$};
    \node (b1) at (0.618,0.3) {$w^1_1$};
    \node (c1) at (1.324,0.3) {$z^0_1$};
    \node (d1) at (-0.662,0.862) {$z^1_1$};
  \node (e1) at (-0.662,-0.262) {$z^2_1$};
   \fill[black] (a) circle (2.5pt);
\fill[black] (b) circle (2.5pt);
\fill[black] (c) circle (2.5pt);
\fill[black] (d) circle (2.5pt);
\fill[black] (e) circle (2.5pt);
\end{tikzpicture}

& \hspace{1cm}
&
\begin{tikzpicture}[->,>=stealth',shorten >=1pt,auto,node distance=2.5cm,main node/.style={circle,draw,font=\sffamily\bfseries \small}, scale=1.35]
    % Draw axes
    \draw [->,thick] (0,-1.2)--(0,1.3) node (yaxis) [above] {$Im(z_2)$}
        |- (-1.75,0)--(1.75,0) node (xaxis) [right] {$Re(z_2)$};

   \node (a) at (-1.618,0){};
    \node (b) at (0.618,0){};
    \node (c) at (1.324,0) {};
    \node (d) at (-0.662,0.562){} ;
  \node (e) at (-0.662,-0.562) {};
  \node (a1) at (-1.618,0.3) {$w^1_2$};
    \node (b1) at (0.618,0.3) {$w^2_2$};
    \node (c1) at (1.324,0.3)  {$z^0_2$};
    \node (d1) at (-0.662,0.862) {$z^1_2$};
  \node (e1) at (-0.662,-0.262) {$z^2_2$};
   \fill[black] (a) circle (2.5pt);
\fill[black] (b) circle (2.5pt);
\fill[black] (c) circle (2.5pt);
\fill[black] (d) circle (2.5pt);
\fill[black] (e) circle (2.5pt);
\end{tikzpicture}
\end{array}$$ Consider the one-parametric family of elements $$g^t(z_1,z_2)=z_1+z_2+\frac{1}{z_1z_2}+\frac{1}{z_1}+\frac{1-t}{z_2} \in L(\Delta^{\circ})$$
for $t \in [0,1)$. Note that $g^0(z_1,z_2)=f_X(z_1,z_2)$ and $g^t(z_1,z_2) \rightarrow f_{Bl_1(\mathbb{P}^2)}(z_1,z_2)$ for $ t \rightarrow 1^-$. Express $Crit(g^t) = \left \{ z^0(t),z^1(t),z^2(t),w^1(t),w^2(t) \right \} \subset (\mathbb{C}^{\ast})^2$. Direct computation gives $$ \begin{array}{cccccc} z^0(t) \rightarrow z^0_{Bl_1(\mathbb{P}^2)} & ; & w^1(t) \rightarrow w^1_{Bl_1(\mathbb{P}^2)}  & ; & w^2(t) \rightarrow w^2(1)=(0,\infty) & ;\end{array}$$
$$ \begin{array}{cccc} z^1(t) \rightarrow z^1_{Bl_1(\mathbb{P}^2)}  & ; & z^2(t) \rightarrow z^2_{Bl_1(\mathbb{P}^2)} & ; \end{array}$$ when $ t \rightarrow 1$. In view of this define
the exceptional map $L : Crit(f_X) \rightarrow Pic(X)$ to be $$ \begin{array}{cccccc} L(z^0)=0 & ; & L(z^1)=H & ; & L(z^2)=2H-E_1-E_2 & ; \end{array}$$
$$ \begin{array}{ccccc} & L(w^1)=H-E_1 & ; & L(w^2)=H-E_2 & ; \end{array}$$

\bigskip

\hspace{-0.6cm} (v) For $\underline{ X=Bl_3(\mathbb{ P}^2)}$ the LG-potential is given by $f_X(z_1,z_2) = z_1+z_2+\frac{1}{z_1z_2}+\frac{1}{z_1}+\frac{1}{z_2}+z_1 z_2$ and the Landau-Ginzburg system is $$ \begin{array}{ccc} z_1 \frac{\partial}{\partial z_1} f_X( z_1,z_2) =
z_1 - \frac{1}{z_1 z_2}-\frac{1}{z_1}+z_1z_2 =0 & ; &  z_2 \frac{\partial}{\partial z_2} f_X( z_1,z_2) =
z_2 - \frac{1}{z_1 z_2} -\frac{1}{z_2}+z_1z_2=0 \end{array} $$ The solution set $Crit(f_X) \subset (\mathbb{C}^{\ast})^2$ is given by $$ \begin{array}{ccccccc} z^0=(1,1) & ; & w^1=(1,-1) & ; & w^2=(-1,1) & ; & w^3=(-1,1) \end{array}$$
$$ \begin{array}{cccc}  z^1=(e^{\frac{2 \pi i}{3}},e^{\frac{2 \pi i}{3}}) & ; & z^2=(e^{\frac{4 \pi i}{3}},e^{\frac{4 \pi i}{3}}) & ; \end{array}$$ Consider the one-parametric family of elements $g^t(z_1,z_2)=
z_1+z_2+\frac{1}{z_1z_2}+\frac{1}{z_1}+\frac{1}{z_2}+(1-t) z_1 z_2$ for $t \in [0,1)$. Note that $g^0(z_1,z_2)=f_X(z_1,z_2)$ and $g^t(z_1,z_2) \rightarrow f_{Bl_2(\mathbb{P}^2)}(z_1,z_2)$ for $ t \rightarrow 1^-$. Express $Crit(g^t) = \left \{ z^0(t),z^1(t),z^2(t),w^1(t),w^2(t),
w^3(t) \right \} \subset (\mathbb{C}^{\ast})^2$. Direct computation gives $$ \begin{array}{cccccc} z^0(t) \rightarrow z^0_{Bl_2(\mathbb{P}^2)} & ; & z^1(t) \rightarrow z^1_{Bl_2(\mathbb{P}^2)}  & ; & z^2(t) \rightarrow z^2_{Bl_2(\mathbb{P}^2)} & ;\end{array}$$ $$ \begin{array}{cccccc}  w^1(t) \rightarrow w^1_{Bl_2(\mathbb{P}^2)}  & ; & w^2(t) \rightarrow w^1_{Bl_2(\mathbb{P}^2)}  & ; & w^3(t) \rightarrow w^3(1)=(0,\infty) & ; \end{array}$$ when $ t \rightarrow 1$. In view of this define
the exceptional map $L : Crit(f_X) \rightarrow Pic(X)$ to be $$ \begin{array}{cccccc} L(z^0)=0 & ; & L(z^1)=H & ; & L(z^2)=2H-E_1-E_2 & ; \end{array}$$
$$ \begin{array}{ccccccc} & L(w^1)=H-E_1 & ; & L(w^2)=H-E_2 & ; & L(w^3)=H-E_3 & ; \end{array}$$

\bigskip

\hspace{-0.6cm} 5.1.2 \bf The Fano $\mathbb{P}^1$-bundles over $\mathbb{P}^2$ case \rm

\bigskip

\hspace{-0.6cm} For  $X=\mathbb{P}(\mathcal{O}_{\mathbb{P}^2} \oplus \mathcal{O}_{\mathbb{P}^2}(k))$ with $k=0,1,2$ the LG-potential is given by
$$ f_X(z_1,z_2,z_3)= z_1+z_2+z_3 + \frac{z_3^k}{z_1z_2} +\frac{1}{z_3} $$ and the Landau-Ginzburg system is $$ \begin{array}{cccc} z_1 \frac{\partial}{\partial z_1} f_X( z_1,z_2) =
z_1 - \frac{z_3^k}{z_1 z_2} =0 & ; &  z_2 \frac{\partial}{\partial z_2} f_X( z_1,z_2) =
z_2 - \frac{z_3^k}{z_1 z_2} =0 & ; \end{array}$$ $$ \begin{array}{cc} z_3 \frac{\partial}{\partial z_3} f_X( z_1,z_2) =
z_3 + \frac{kz_3^k}{z_1 z_2}-\frac{1}{z_3} =0 & ; \end{array} $$

\bigskip

\hspace{-0.6cm} (i) For $k=0$ the solution set $Crit(f_X) \subset (\mathbb{C}^{\ast})^3$ is given by $$ \begin{array}{cccc} z^{00} = (1,1,1) & ; & z^{01}\approx (1,1,-1) & ; \end{array}$$ $$ \begin{array}{ccccc}
z^{10} = (e^{\frac{2 \pi i}{3}}, e^{\frac{2 \pi i}{3}},1) & ;  z^{11} = ( e^{\frac{2 \pi i}{3}}, e^{\frac{2 \pi i}{3}},-1) & ; \end{array}$$
$$ \begin{array}{ccccc} z^{20} =(e^{\frac{4 \pi i}{3}}, e^{\frac{4 \pi i}{3}},1) & ; & z^{21}=(e^{\frac{4 \pi i}{3}}, e^{\frac{4 \pi i}{3}},-1)   & ; \end{array}$$

\bigskip

\hspace{-0.6cm} (ii) For $k=1$ the solution set $Crit(f_X) \subset (\mathbb{C}^{\ast})^3$ is given by
$$ \begin{array}{ccccc} z^{00} \approx (0.86,0.86,0.65) & ; & z^{11} \approx (-0.86,-0.86,-0.65)  & ; \end{array}$$
$$ \begin{array}{ccccc} z^{10} \approx (-0.36+i,-0.36+i,1.07-0.6i) & ; & z^{21}\approx (0.36-i,0.36-i,-1.07+0.6i) & ; \end{array}$$
$$ \begin{array}{ccccc} z^{20}  \approx (-0.36-i,-0.36-i,1.07+0.6i) & ;  & z^{01}\approx (0.36+i,0.36+i,-1.07-0.6i) & ; \end{array}$$ The following is a graphical description of the $z_3$-plane:
$$ 
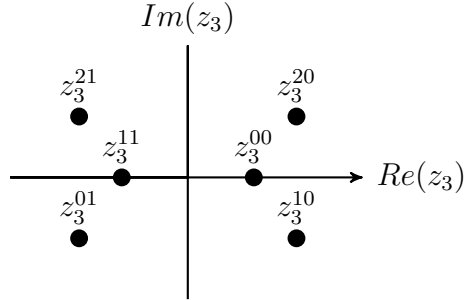
\begin{tikzpicture}[->,>=stealth',shorten >=1pt,auto,node distance=2.5cm,main node/.style={circle,draw,font=\sffamily\bfseries \small}, scale=1.35]
    % Draw axes
    \draw [->,thick] (0,-1.2)--(0,1.3) node (yaxis) [above] {$Im(z_3)$}
        |- (-1.75,0)--(1.75,0) node (xaxis) [right] {$Re(z_3)$};

   \node (a) at (0.65,0){};
    \node (b) at (-0.65,0){};
    \node (c) at (1.07,0.6) {};
    \node (d) at (-1.07,0.6) {};
  \node (e) at (1.07,-0.6) {};
 \node (f) at (-1.07,-0.6) {};
  \node (a1) at  (0.65,0.3) {$z^{00}_3$};
    \node (b1) at  (-0.65,0.3) {$z^{11}_3$};
    \node (c1) at (1.07,0.9) {$z^{20}_3$};
    \node (d1) at  (-1.07,0.9) {$z^{21}_3$};
  \node (e1) at (1.07,-0.3) {$z^{10}_3$};
 \node (f1) at (-1.07,-0.3)  {$z^{01}_3$};
   \fill[black] (a) circle (2.5pt);
\fill[black] (b) circle (2.5pt);
\fill[black] (c) circle (2.5pt);
\fill[black] (d) circle (2.5pt);
\fill[black] (e) circle (2.5pt);
\fill[black] (f) circle (2.5pt);
\end{tikzpicture}$$

\hspace{-0.6cm} (iii) For $k=2$ the solution set $Crit(f_X) \subset (\mathbb{C}^{\ast})^3$ is given by
$$ \begin{array}{ccccc} z^{00} \approx (0.66,0.66,0.53) & ; & z^{21}  \approx (4.11,4.11,-8.35)  & ; \end{array}$$
$$ \begin{array}{cccc} z^{10} \approx (0.11-0.85i,0.11-0.85i,0.43+0.67i) & ; & z^{11} \approx  (-0.5-0.4i,-0.5-0.4i,-0.53+0.24i)& ; \end{array}$$
 $$ \begin{array}{ccccc} z^{20} \approx (0.11+0.85i,0.11+0.85i,0.43-0.67i) & ;  z^{01}  \approx (-0.5+0.4i,-0.5+0.4i,-0.53-0.24i)  & ; \end{array}$$
The following is a graphical description of the $z_3$-plane:
$$ 
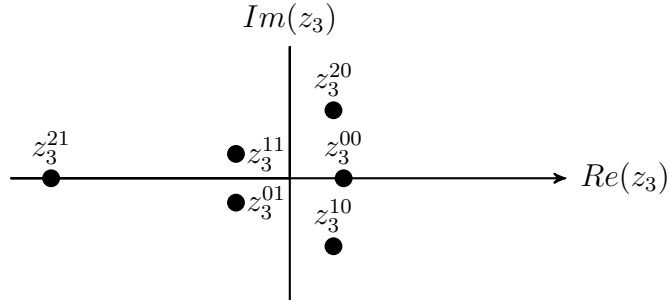
\begin{tikzpicture}[->,>=stealth',shorten >=1pt,auto,node distance=2.5cm,main node/.style={circle,draw,font=\sffamily\bfseries \small}, scale=1.35]
    % Draw axes
    \draw [->,thick] (0,-1.2)--(0,1.3) node (yaxis) [above] {$Im(z_3)$}
        |- (-2.75,0)--(2.75,0) node (xaxis) [right] {$Re(z_3)$};

   \node (a) at (0.53,0){};
    \node (b) at (-2.35,0){};
    \node (c) at (0.43,0.67){} ;
    \node (d) at (-0.53,0.24){} ;
  \node (e) at (0.43,-0.67) {};
 \node (f) at (-0.53,-0.24) {};
  \node (a1) at  (0.53,0.3) {$z^{00}_3$};
    \node (b1) at  (-2.35,0.3) {$z^{21}_3$};
    \node (c1) at (0.43,0.97) {$z^{20}_3$};
    \node (d1) at  (-0.23,0.24) {$z^{11}_3$};
  \node (e1) at (0.43,-0.37) {$z^{10}_3$};
 \node (f1) at (-0.23,-0.24)  {$z^{01}_3$};
   \fill[black] (a) circle (2.5pt);
\fill[black] (b) circle (2.5pt);
\fill[black] (c) circle (2.5pt);
\fill[black] (d) circle (2.5pt);
\fill[black] (e) circle (2.5pt);
\fill[black] (f) circle (2.5pt);
\end{tikzpicture}$$ In the three cases we define $L : Crit(f_X) \rightarrow Pic(X)$ by $L(z^{lm})=l \pi^{\ast}+m\xi$ where $0 \leq l \leq 2 $ and $0 \leq m \leq 1$. This is justified as follows: consider the
one-parametric family of elements $$ g^t(z_1,z_2,z_3) = (1-t)z_1+(1-t)z_2 +z_3 +\frac{(1-t)z^k}{z_1 z_2} + \frac{1}{z_3} \in L(\Delta^{\circ})$$ for $t \in [0,1)$.
Note that $g^0(z)=f_X(z)$ and $g^t(z) \rightarrow f_{\mathbb{P}^1}(z)$ for $ t \rightarrow 1^-$. Express $$Crit(g^t) = \left \{ z^{ij}(t) \vert 0 \leq i \leq 2,0 \leq j \leq 1 \right \} \subset (\mathbb{C}^{\ast})^3$$ Direct computation shows
$$ \begin{array}{cccc} z^{00}(t),z^{10}(t),z^{20}(t) \rightarrow \widetilde{z}^0_{\mathbb{P}^1}=(0,0,1)  & ; & z^{01}(t),z^{11}(t),z^{21}(t) \rightarrow \widetilde{z}^1_{\mathbb{P}^1}=(0,0,-1) & ; \end{array}$$ when $t \rightarrow 1$.
On the other hand consider the one-parametric family of elements $$ h^t(z_1,z_2,z_3) = z_1+z_2 +z_3 +\frac{ e^{2 \pi i t} z^k}{z_1 z_2}  + \frac{1}{z_3}$$ for $t \in [0,1)$. Computation shows
$$ \begin{array}{cccccc} z^{00}(t) \rightarrow z^{10}  & ; & z^{10}(t) \rightarrow z^{20} & ; & z^{20}(t) \rightarrow z^{00} ; \end{array}$$
$$ \begin{array}{cccccc} z^{01}(t) \rightarrow z^{11}  & ; & z^{11}(t) \rightarrow z^{21} & ; & z^{21}(t) \rightarrow z^{01} ; \end{array}$$

\bigskip

\hspace{-0.6cm} 5.1.3 \bf The Fano $\mathbb{P}^2$-bundles over $\mathbb{P}^1$ case \rm

\bigskip

\hspace{-0.6cm} For  $X=\mathbb{P}(\mathcal{O}_{\mathbb{P}^1} \oplus \mathcal{O}_{\mathbb{P}^1} \oplus \mathcal{O}_{\mathbb{P}^1}(k))$ with $k=0,1$ the LG-potential is given by
$$ f_X(z_1,z_2,z_3)= z_1+z_2+z_3 + \frac{1}{z_1z_2} +\frac{1}{z_1^k z_2^k z_3} $$ For $k=0$ this is the same as Example 4.2.(i). For $k=1$ the Landau-Ginzburg system is $$ \begin{array}{cccc} z_1 \frac{\partial}{\partial z_1} f_X( z_1,z_2) =
z_1 - \frac{1}{z_1 z_2} -\frac{1}{z_1z_2z_3}=0 & ; &  z_2 \frac{\partial}{\partial z_2} f_X( z_1,z_2) =
z_2 - \frac{1}{z_1 z_2}  -\frac{1}{z_1z_2z_3}=0 & ; \end{array}$$ $$ \begin{array}{cc} z_3 \frac{\partial}{\partial z_3} f_X( z_1,z_2) =
z_3  -\frac{1}{z_1z_2}-\frac{1}{z_1^k z_2^k z_3} =0 & ; \end{array} $$ and the solution set $Crit(f_X) \subset (\mathbb{C}^{\ast})^2$ is given by

$$ \begin{array}{ccccc} z^{00} \approx (1.3,1.3,0.7) & ; & z^{10}  \approx (0.6,0.6,-1.4)  & ; \end{array}$$
$$ \begin{array}{cccc} z^{01} \approx (-0.3+1.1i,-0.3+1.1i,-0.2+0.7i) & ; & z^{11} \approx  (-0.6+0.5i,-0.6+0.5i,-0.8-0.7i) & ; \end{array}$$
 $$ \begin{array}{ccccc} z^{02} \approx (-0.6-0.5i,-0.6-0.5i,-0.8-0.7i) & ;  z^{12}  \approx(-0.3-1.1i,-0.3-1.1i,-0.2-0.7i)  & ; \end{array}
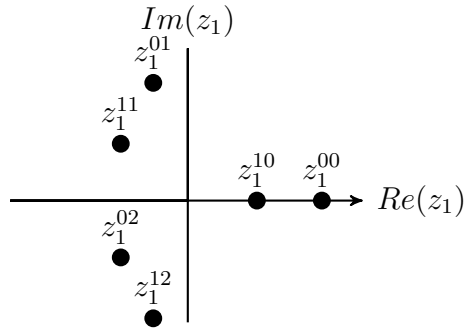
$$

$$ \begin{tikzpicture}[->,>=stealth',shorten >=1pt,auto,node distance=2.5cm,main node/.style={circle,draw,font=\sffamily\bfseries \small}, scale=1.35]
    % Draw axes
    \draw [->,thick] (0,-1.2)--(0,1.5) node (yaxis) [above] {$Im(z_1)$}
        |- (-1.75,0)--(1.75,0) node (xaxis) [right] {$Re(z_1)$};

    \node (a) at (1.32,0){};
    \node (b) at (-0.34,1.16){};
    \node (c) at (-0.66,-0.56) {};
  \node (a1) at (1.32,0.3) {$z^{00}_1$};
    \node (b1) at (-0.34,1.46) {$z^{01}_1$};
    \node (c1) at (-0.66,-0.26) {$z^{02}_1$};
 \node (a2) at (0.68,0){};
    \node (b2) at (-0.66,0.56){};
    \node (c2) at (-0.34,-1.16) {};
  \node (a3) at (0.68,0.3) {$z^{10}_1$};
    \node (b3) at (-0.66,0.86) {$z^{11}_1$};
    \node (c3) at (-0.34,-0.86) {$z^{12}_1$};
   \fill[black] (a) circle (2.5pt);
\fill[black] (b) circle (2.5pt);
\fill[black] (c) circle (2.5pt);
  \fill[black] (a2) circle (2.5pt);
\fill[black] (b2) circle (2.5pt);
\fill[black] (c2) circle (2.5pt);
\end{tikzpicture} $$ By similar considerations to those of Example 4.2 we define $L : Crit(f_X) \rightarrow Pic(X)$ to be $L(z^{lm})=l\pi^{\ast}H+m\xi$ where $0 \leq l \leq 1 $ and $0 \leq m \leq 2$.

\bigskip

\hspace{-0.6cm}  5.1.4 \bf The Fano $\mathbb{P}^1$-bundles over $\mathbb{P}^1 \times \mathbb{P}^1$ case \rm

\bigskip

\hspace{-0.6cm} For $X=\mathbb{P}(\mathcal{O}_{\mathbb{P}^1 \times \mathbb{P}^1} \oplus \mathcal{O}_{\mathbb{P}^1 \times \mathbb{P}^1}(k_1,k_2))$ with $(k_1,k_2)=(0,0),(1,1),(1,-1)$ the LG-potential is given by
$$ f_X(z_1,z_2,z_3)= z_1+z_2+z_3 + \frac{z_1^{k_1}}{z_2} +\frac{z_1^{k_2}}{z_3}+\frac{1}{z_1} $$ and the Landau-Ginzburg system is $$ \begin{array}{cccc} z_1 \frac{\partial}{\partial z_1} f_X( z_1,z_2) =
z_1 +k_1 \frac{z_1^{k_1}}{ z_2} +k_2 \frac{z_1^{k_2}}{z_3} -\frac{1}{z_1} =0 & ; &  z_2 \frac{\partial}{\partial z_2} f_X( z_1,z_2) =
z_2 - \frac{z_1^{k_1}}{z_2} =0 & ; \end{array}$$ $$ \begin{array}{cc} z_3 \frac{\partial}{\partial z_3} f_X( z_1,z_2) =
z_3 -\frac{z_1^{k_2}}{z_3} =0 & ; \end{array} $$

\hspace{-0.6cm} (i) For $(k_1,k_2)=(0,0)$ the solution set $Crit(f_X) \subset (\mathbb{C}^{\ast})^2$ is given by
$$ \begin{array}{cccc} z^{000} = (1,1,1) & ; & z^{001}= (1,1,-1) & ; \end{array}$$
$$ \begin{array}{cccc} z^{010} = (1,-1,1) & ; & z^{011}= (1,-1,-1) & ; \end{array}$$
$$ \begin{array}{cccc} z^{100} = (-1,1,1) & ; & z^{101}= (-1,1,-1) & ; \end{array}$$
$$ \begin{array}{cccc} z^{110} = (-1,-1,1) & ; & z^{111}= (-1,-1,-1) & ; \end{array}$$

\hspace{-0.6cm} (ii) For $(k_1,k_2)=(1,1)$ the solution set $Crit(f_X) \subset (\mathbb{C}^{\ast})^2$ is given by
$$ \begin{array}{cccc} z^{000} \approx  (0.51,0.71,0.71) & ; & z^{001} \approx  (-0.47,-0.3+0.75i,-0.3+0.75i) & ; \end{array}$$
$$ \begin{array}{cccc} z^{010} \approx (1,-1,1) & ; & z^{011} \approx (-1,i,-i) & ; \end{array}$$
$$ \begin{array}{cccc} z^{100} \approx (1,1,-1) & ; & z^{101} \approx (-1,-i,i)  & ; \end{array}$$
$$ \begin{array}{cccc} z^{110} \approx (4.43,-2.1,-2.1) & ; & z^{111} \approx (-0.47,-0.3-0.75i,-0.3-0.75i)& ; \end{array}
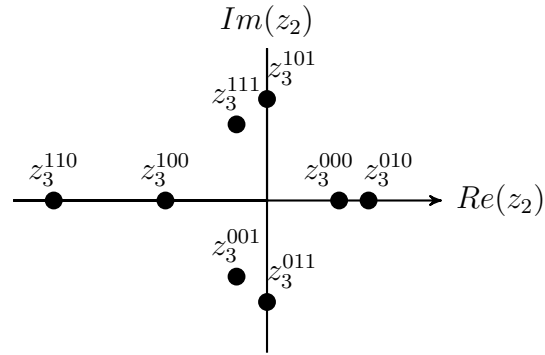
$$
$$ \begin{tikzpicture}[->,>=stealth',shorten >=1pt,auto,node distance=2.5cm,main node/.style={circle,draw,font=\sffamily\bfseries \small}, scale=1.35]
    % Draw axes
    \draw [->,thick] (0,-1.5)--(0,1.5) node (yaxis) [above] {$Im(z_2)$}
        |- (-2.5,0)--(1.75,0) node (xaxis) [right] {$Re(z_2)$};

    \node (a) at (0.71,0){};
    \node (b) at (1,0){};
    \node (c) at (-1,0) {};
\node (d) at (-2.1,0) {};
  \node (a1) at (0.61,0.3) {$z^{000}_3$};
    \node (b1) at (1.2,0.3) {$z^{010}_3$};
    \node (c1) at (-1,0.3) {$z^{100}_3$};
    \node (d1) at (-2.1,0.3) {$z^{110}_3$};
 \node (a2) at (-0.3,0.75){};
    \node (b2) at (0,1){};
    \node (c2) at (0,-1) {};
 \node (d2) at (-0.3,-0.75){} ;
  \node (a3) at (-0.3,1.05) {$z^{111}_3$};
    \node (b3) at (0.25,1.3) {$z^{101}_3$};
    \node (c3) at (0.25,-0.7) {$z^{011}_3$};
    \node (d3) at (-0.3,-0.45) {$z^{001}_3$};
   \fill[black] (a) circle (2.5pt);
\fill[black] (b) circle (2.5pt);
\fill[black] (c) circle (2.5pt);
\fill[black] (d) circle (2.5pt);
  \fill[black] (a2) circle (2.5pt);
\fill[black] (b2) circle (2.5pt);
\fill[black] (c2) circle (2.5pt);
\fill[black] (d2) circle (2.5pt);
\end{tikzpicture} $$

\hspace{-0.6cm} (iii) For $(k_1,k_2)=(1,-1)$ the solution set $Crit(f_X) \subset (\mathbb{C}^{\ast})^2$ is given by
$$ \begin{array}{cccc} z^{000} \approx  (1,1,1) & ; & z^{101} \approx  (-1,i,i) & ; \end{array}$$
$$ \begin{array}{cccc} z^{010} \approx (0.38,0.61,-1.61) & ; & z^{111} \approx (-0.5-0.866i,-0.5+0.866i,-0.5-0.866i) & ; \end{array}$$
$$ \begin{array}{cccc} z^{100} \approx (2.61,-1.61,0.61) & ; & z^{001} \approx (-0.5+0.866i,-0.5-0.886i,-0.5+0.866i)  & ; \end{array}$$
$$ \begin{array}{cccc} z^{110} \approx (1,-1,-1) & ; & z^{011} \approx (-1,-i,-i)& ; \end{array}
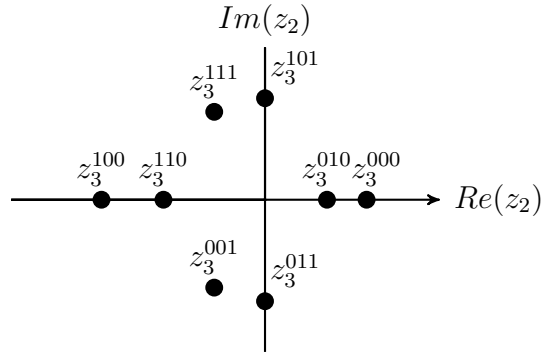
$$
$$ \begin{tikzpicture}[->,>=stealth',shorten >=1pt,auto,node distance=2.5cm,main node/.style={circle,draw,font=\sffamily\bfseries \small}, scale=1.35]
    % Draw axes
    \draw [->,thick] (0,-1.5)--(0,1.5) node (yaxis) [above] {$Im(z_2)$}
        |- (-2.5,0)--(1.75,0) node (xaxis) [right] {$Re(z_2)$};

    \node (a) at (1,0){};
    \node (b) at (-1.61,0){};
    \node (c) at (0.61,0) {};
\node (d) at (-1,0) {};
  \node (a1) at (1.1,0.3) {$z^{000}_3$};
    \node (b1) at (-1.61,0.3) {$z^{100}_3$};
    \node (c1) at (0.61,0.3) {$z^{010}_3$};
    \node (d1) at (-1,0.3) {$z^{110}_3$};
 \node (a2) at (0,1){};
    \node (b2) at (-0.5,-0.866){};
    \node (c2) at (-0.5,0.866){};
 \node (d2) at (0,-1) {};
  \node (a3) at (0.3,1.3) {$z^{101}_3$};
    \node (b3) at (-0.5,-0.566) {$z^{001}_3$};
    \node (c3) at (-0.5,1.166) {$z^{111}_3$};
    \node (d3) at (0.3,-0.7) {$z^{011}_3$};
   \fill[black] (a) circle (2.5pt);
\fill[black] (b) circle (2.5pt);
\fill[black] (c) circle (2.5pt);
\fill[black] (d) circle (2.5pt);
  \fill[black] (a2) circle (2.5pt);
\fill[black] (b2) circle (2.5pt);
\fill[black] (c2) circle (2.5pt);
\fill[black] (d2) circle (2.5pt);
\end{tikzpicture} $$ We define $L : Crit(f_X) \rightarrow Pic(X)$ by $L(z^{lmn})=l\pi^{\ast}H_1+m\pi^{\ast}H_2+n\xi$ for $0 \leq l,m,n \leq 1 $.

\subsection{The $M$-aligned property:}
\label{s:EMD} In this sub-section we verify the $M$-aligned property for the exceptional maps defined in sub-section 5.1. Note that, by definition, verifying the $M$-aligned property requires to show that
$\widetilde{Q}^D(\mathcal{E}) \subset \widetilde{Q}(D)$ for any $D \in Div(\mathcal{E})$. We thus compute $\widetilde{Q}(D)$ for any $D \in Div(\mathcal{E})$ and compare it to the quivers
described in Example 3.11-14. Let us note that the graphs of the quivers $ \widetilde{Q}(D)$ given below do not represent the actual path curve of the monodromy, but serve as a schematic description of the underlying 
automorphisms. The bold arrows represent the edges of the sub-quiver $\widetilde{Q}^D(\mathcal{E}) \subset \widetilde{Q}(D)$.

\bigskip

\hspace{-0.6cm} 5.2.1 \bf The Del-Pezzo surface case \rm

\bigskip

\hspace{-0.6cm} (i) For $\underline{ X=\mathbb{ P}^2}$ the monodromy quivers are:

$$
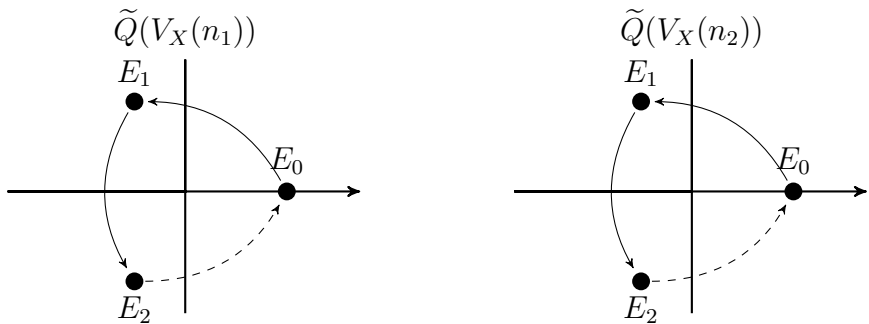
\begin{array}{ccc} \begin{tikzpicture}[->,>=stealth',shorten >=1pt,auto,node distance=2.5cm,main node/.style={circle,draw,font=\sffamily\bfseries \small}, scale=1.35]
    % Draw axes
    \draw [->,thick] (0,-1.2)--(0,1.3) node (yaxis) [above] {$\widetilde{Q}(V_X(n_1))$}
        |- (-1.75,0)--(1.75,0) node (xaxis) [right] {};

    \node (a) at (1,0){};
    \node (b) at (-0.5,0.886){};
    \node (c) at (-0.5,-0.886){} ;
  \node (a1) at (1,0.3) {$E_0$};
    \node (b1) at (-0.5,1.166) {$E_1$};
    \node (c1) at (-0.5,-1.166) {$E_2$};
   \fill[black] (a) circle (2.5pt);
\fill[black] (b) circle (2.5pt);
\fill[black] (c) circle (2.5pt);

 \path[every node/.style={font=\sffamily\small}]

    (a) edge [bend right] (b)
     (b) edge [bend right] (c)
(c) edge [dashed,->,>=stealth',shorten >=1pt,auto,bend right]  (a);

\end{tikzpicture}

& \hspace{1cm}
&
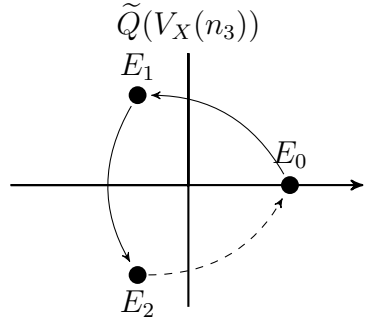
\begin{tikzpicture}[->,>=stealth',shorten >=1pt,auto,node distance=2.5cm,main node/.style={circle,draw,font=\sffamily\bfseries \small}, scale=1.35]
    % Draw axes
    \draw [->,thick] (0,-1.2)--(0,1.3) node (yaxis) [above] {$\widetilde{Q}(V_X(n_2))$}
        |- (-1.75,0)--(1.75,0) node (xaxis) [right] {};

    \node (a) at (1,0){};
    \node (b) at (-0.5,0.886){};
    \node (c) at (-0.5,-0.886){} ;
 \node (a1) at (1,0.3) {$E_0$};
    \node (b1) at (-0.5,1.166) {$E_1$};
    \node (c1) at (-0.5,-1.166) {$E_2$};
   \fill[black] (a) circle (2.5pt);
\fill[black] (b) circle (2.5pt);
\fill[black] (c) circle (2.5pt);

 \path[every node/.style={font=\sffamily\small}]

    (a) edge [bend right] (b)
     (b) edge [bend right] (c)
(c) edge  [dashed,->,>=stealth',shorten >=1pt,auto,bend right]  (a);

\end{tikzpicture}

\end{array}%
$$

\bigskip

\hspace{-0.6cm} (ii) For $\underline{ X=\mathbb{ P}^1 \times \mathbb{P}^1}$ the monodromy quivers are:

\bigskip

$$ \begin{array}{ccc} 
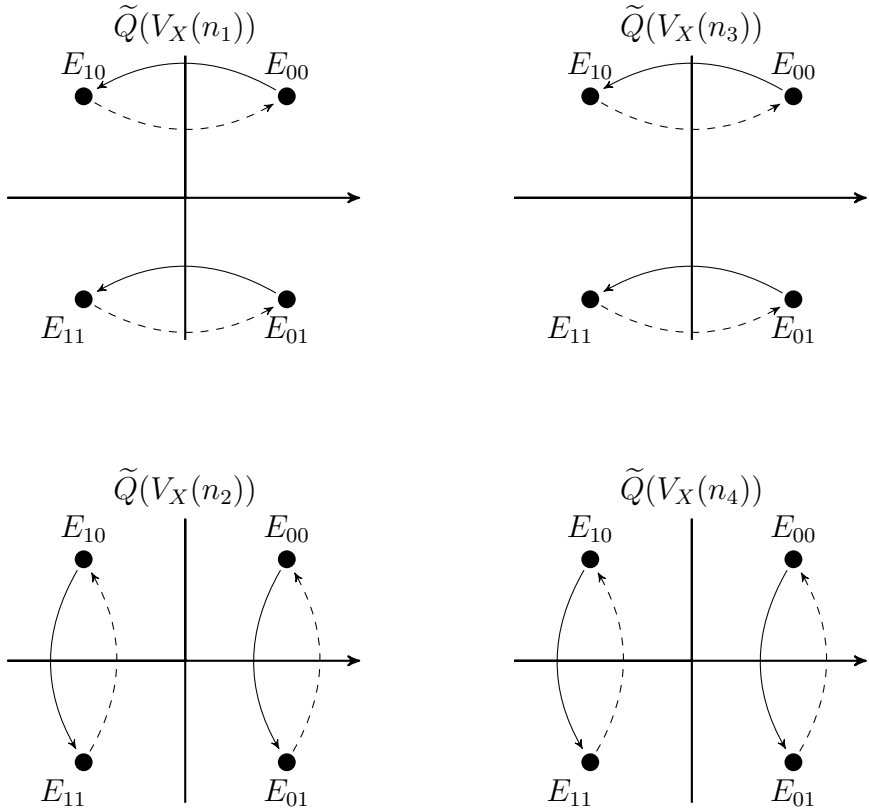
\begin{tikzpicture}[->,>=stealth',shorten >=1pt,auto,node distance=2.5cm,main node/.style={circle,draw,font=\sffamily\bfseries \small}, scale=1.35]

    % Draw axes
    \draw [->,thick] (0,-1.4)--(0,1.4) node (yaxis) [above] {$\widetilde{Q}(V_X(n_1))$}
        |- (-1.75,0)--(1.75,0) node (xaxis) [right] {};

    \node (a) at (1,1){};
    \node (b) at (-1,1){};
    \node (c) at (1,-1) {};
\node (d) at (-1,-1){};
  \node (a1) at (1,1.3) {$E_{00}$};
    \node (b1) at (-1,1.3) {$E_{10}$};
    \node (c1) at (1,-1.3) {$E_{01}$};
\node (d1) at (-1.2,-1.3) {$E_{11}$};
   \fill[black] (a) circle (2.5pt);
\fill[black] (b) circle (2.5pt);
\fill[black] (c) circle (2.5pt);
\fill[black] (d) circle (2.5pt);

 \path[every node/.style={font=\sffamily\small}]

    (a) edge [bend right]  (b)
     (b) edge  [dashed,->,>=stealth',shorten >=1pt,auto,bend right]   (a)
(c) edge [bend right] (d)
(d) edge [dashed,->,>=stealth',shorten >=1pt,auto,bend right]   (c);

\end{tikzpicture}
& \hspace{1cm}
&
 \begin{tikzpicture}[->,>=stealth',shorten >=1pt,auto,node distance=2.5cm,main node/.style={circle,draw,font=\sffamily\bfseries \small}, scale=1.35]

    % Draw axes
    \draw [->,thick] (0,-1.4)--(0,1.4) node (yaxis) [above] {$\widetilde{Q}(V_X(n_3))$}
        |- (-1.75,0)--(1.75,0) node (xaxis) [right] {};

    \node (a) at (1,1){};
    \node (b) at (-1,1){};
    \node (c) at (1,-1) {};
\node (d) at (-1,-1){};
 \node (a1) at (1,1.3) {$E_{00}$};
    \node (b1) at (-1,1.3) {$E_{10}$};
    \node (c1) at (1,-1.3) {$E_{01}$};
\node (d1) at (-1.2,-1.3) {$E_{11}$};
   \fill[black] (a) circle (2.5pt);
\fill[black] (b) circle (2.5pt);
\fill[black] (c) circle (2.5pt);
\fill[black] (d) circle (2.5pt);

 \path[every node/.style={font=\sffamily\small}]

    (a) edge [bend right]  (b)
     (b) edge  [dashed,->,>=stealth',shorten >=1pt,auto,bend right]   (a)
(c) edge [bend right]  (d)
(d) edge  [dashed,->,>=stealth',shorten >=1pt,auto,bend right]   (c);

\end{tikzpicture}
\end{array}
$$

$$ \begin{array}{ccc} \begin{tikzpicture}[->,>=stealth',shorten >=1pt,auto,node distance=2.5cm,main node/.style={circle,draw,font=\sffamily\bfseries \small}, scale=1.35]

    % Draw axes
    \draw [->,thick] (0,-1.4)--(0,1.4) node (yaxis) [above] {$\widetilde{Q}(V_X(n_2))$}
        |- (-1.75,0)--(1.75,0) node (xaxis) [right] {};

    \node (a) at (1,1){};
    \node (b) at (-1,1){};
    \node (c) at (1,-1){} ;
\node (d) at (-1,-1){};
 \node (a1) at (1,1.3) {$E_{00}$};
    \node (b1) at (-1,1.3) {$E_{10}$};
    \node (c1) at (1,-1.3) {$E_{01}$};
\node (d1) at (-1.2,-1.3) {$E_{11}$};
   \fill[black] (a) circle (2.5pt);
\fill[black] (b) circle (2.5pt);
\fill[black] (c) circle (2.5pt);
\fill[black] (d) circle (2.5pt);

\path[every node/.style={font=\sffamily\small}]

    (a) edge [bend right] (c)
     (b) edge [bend right]  (d)
(c) edge  [dashed,->,>=stealth',shorten >=1pt,auto,bend right] (a)
(d) edge  [dashed,->,>=stealth',shorten >=1pt,auto,bend right]  (b);

\end{tikzpicture}
& \hspace{1cm}
&
\begin{tikzpicture}[->,>=stealth',shorten >=1pt,auto,node distance=2.5cm,main node/.style={circle,draw,font=\sffamily\bfseries \small}, scale=1.35]

    % Draw axes
    \draw [->,thick] (0,-1.4)--(0,1.4) node (yaxis) [above] {$\widetilde{Q}(V_X(n_4))$}
        |- (-1.75,0)--(1.75,0) node (xaxis) [right] {};

    \node (a) at (1,1){};
    \node (b) at (-1,1){};
    \node (c) at (1,-1){} ;
\node (d) at (-1,-1){};
 \node (a1) at (1,1.3) {$E_{00}$};
    \node (b1) at (-1,1.3) {$E_{10}$};
    \node (c1) at (1,-1.3) {$E_{01}$};
\node (d1) at (-1.2,-1.3) {$E_{11}$};
   \fill[black] (a) circle (2.5pt);
\fill[black] (b) circle (2.5pt);
\fill[black] (c) circle (2.5pt);
\fill[black] (d) circle (2.5pt);

 \path[every node/.style={font=\sffamily\small}]
(a) edge [bend right]  (c)
     (b) edge [bend right]  (d)
  (c) edge  [dashed,->,>=stealth',shorten >=1pt,auto,bend right]   (a)
(d) edge  [dashed,->,>=stealth',shorten >=1pt,auto,bend right]  (b);

\end{tikzpicture}
\end{array}
$$

\bigskip

\hspace{-0.6cm} (iii) For $\underline{ X=Bl_1(\mathbb{ P}^2)}$ the monodromy quivers are:
$$ \begin{array}{cccc}
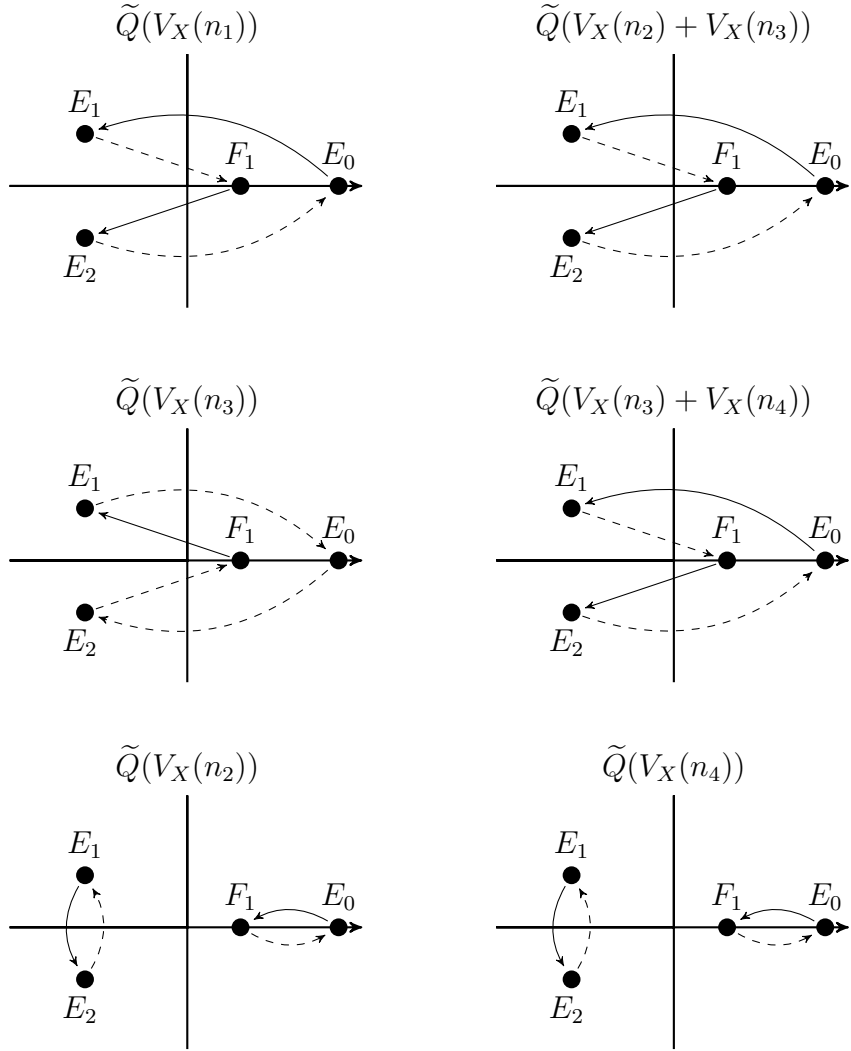
\begin{tikzpicture}[->,>=stealth',shorten >=1pt,auto,node distance=2.5cm,main node/.style={circle,draw,font=\sffamily\bfseries \small}, scale=1.35]
    % Draw axes
    \draw [->,thick] (0,-1.2)--(0,1.3) node (yaxis) [above] {$\widetilde{Q}(V_X(n_1))$}
        |- (-1.75,0)--(1.75,0);

    \node (a) at (0.524,0) {};
    \node (b) at (-1.007,0.513){};
    \node (c) at (-1.007,-0.513){} ;
    \node (d) at (1.490,0){} ;

\node (a1) at (0.524,0.3) {$F_1$};
    \node (b1) at (-1.007,0.813) {$E_1$};
    \node (c1) at (-1.057,-0.813) {$E_2$} ;
    \node (d1) at (1.490,0.3) {$E_0$};

   \fill[black] (a) circle (2.5pt);
\fill[black] (b) circle (2.5pt);
\fill[black] (c) circle (2.5pt);
\fill[black] (d) circle (2.5pt);

 \path[every node/.style={font=\sffamily\small}]

(c) edge [dashed,->,bend right] (d)
      (b) edge [dashed,->] (a)
(d) edge [bend right] (b)
 (a)  edge (c);

\end{tikzpicture}
& \hspace{1cm}
&
\begin{tikzpicture}[->,>=stealth',shorten >=1pt,auto,node distance=2.5cm,main node/.style={circle,draw,font=\sffamily\bfseries \small}, scale=1.35]
    % Draw axes
    \draw [->,thick] (0,-1.2)--(0,1.3) node (yaxis) [above] {$\widetilde{Q}(V_X(n_2)+V_X(n_3))$}
        |- (-1.75,0)--(1.75,0);

    \node (a) at (0.524,0) {};
    \node (b) at (-1.007,0.513){};
    \node (c) at (-1.007,-0.513){} ;
    \node (d) at (1.490,0){} ;

\node (a1) at (0.524,0.3) {$F_1$};
    \node (b1) at (-1.007,0.813) {$E_1$};
    \node (c1) at (-1.057,-0.813) {$E_2$} ;
    \node (d1) at (1.490,0.3) {$E_0$};

   \fill[black] (a) circle (2.5pt);
\fill[black] (b) circle (2.5pt);
\fill[black] (c) circle (2.5pt);
\fill[black] (d) circle (2.5pt);

 \path[every node/.style={font=\sffamily\small}]

(c) edge [dashed,->,bend right] (d)
      (b) edge [dashed,->] (a)
(d) edge [bend right] (b)
 (a)  edge (c);

\end{tikzpicture}
\end{array}$$

$$\begin{array}{ccc}
\begin{tikzpicture}[->,>=stealth',shorten >=1pt,auto,node distance=2.5cm,main node/.style={circle,draw,font=\sffamily\bfseries \small}, scale=1.35]
    % Draw axes
    \draw [->,thick] (0,-1.2)--(0,1.3) node (yaxis) [above] {$\widetilde{Q}(V_X(n_3))$}
        |- (-1.75,0)--(1.75,0);

     \node (a) at (0.524,0) {};
    \node (b) at (-1.007,0.513){};
    \node (c) at (-1.007,-0.513){} ;
    \node (d) at (1.490,0){} ;

\node (a1) at (0.524,0.3) {$F_1$};
    \node (b1) at (-1.007,0.813) {$E_1$};
    \node (c1) at (-1.057,-0.813) {$E_2$} ;
    \node (d1) at (1.490,0.3) {$E_0$};

   \fill[black] (a) circle (2.5pt);
\fill[black] (b) circle (2.5pt);
\fill[black] (c) circle (2.5pt);
\fill[black] (d) circle (2.5pt);

 \path[every node/.style={font=\sffamily\small}]

(d) edge [dashed,->,bend left] (c)
 (b) edge [dashed,->,bend left] (d)
(c) edge [dashed,->]  (a)
  (a) edge  (b);

\end{tikzpicture}
& \hspace{1cm}
&
\begin{tikzpicture}[->,>=stealth',shorten >=1pt,auto,node distance=2.5cm,main node/.style={circle,draw,font=\sffamily\bfseries \small}, scale=1.35]
    % Draw axes
    \draw [->,thick] (0,-1.2)--(0,1.3) node (yaxis) [above] {$\widetilde{Q}(V_X(n_3)+V_X(n_4))$}
        |- (-1.75,0)--(1.75,0);

      \node (a) at (0.524,0) {};
    \node (b) at (-1.007,0.513){};
    \node (c) at (-1.007,-0.513){} ;
    \node (d) at (1.490,0){} ;

\node (a1) at (0.524,0.3) {$F_1$};
    \node (b1) at (-1.007,0.813) {$E_1$};
    \node (c1) at (-1.057,-0.813) {$E_2$} ;
    \node (d1) at (1.490,0.3) {$E_0$};

   \fill[black] (a) circle (2.5pt);
\fill[black] (b) circle (2.5pt);
\fill[black] (c) circle (2.5pt);
\fill[black] (d) circle (2.5pt);

 \path[every node/.style={font=\sffamily\small}]

(d) edge [bend right] (b)
 (a)  edge (c)

(c) edge [dashed,->,bend right] (d)
      (b) edge [dashed,->] (a);

\end{tikzpicture}
\end{array}$$

$$ \begin{array}{ccc}
 \begin{tikzpicture}[->,>=stealth',shorten >=1pt,auto,node distance=2.5cm,main node/.style={circle,draw,font=\sffamily\bfseries \small}, scale=1.35]
    % Draw axes
    \draw [->,thick] (0,-1.2)--(0,1.3) node (yaxis) [above] {$\widetilde{Q}(V_X(n_2))$}
        |- (-1.75,0)--(1.75,0);

      \node (a) at (0.524,0) {};
    \node (b) at (-1.007,0.513){};
    \node (c) at (-1.007,-0.513){} ;
    \node (d) at (1.490,0){} ;

\node (a1) at (0.524,0.3) {$F_1$};
    \node (b1) at (-1.007,0.813) {$E_1$};
    \node (c1) at (-1.057,-0.813) {$E_2$} ;
    \node (d1) at (1.490,0.3) {$E_0$};

   \fill[black] (a) circle (2.5pt);
\fill[black] (b) circle (2.5pt);
\fill[black] (c) circle (2.5pt);
\fill[black] (d) circle (2.5pt);

\path[every node/.style={font=\sffamily\small}]
(a) edge [dashed,->,bend right]  (d)
     (d) edge [bend right] (a)
  (c) edge [dashed,->,bend right] (b)
(b) edge [bend right] (c);

\end{tikzpicture}
 & \hspace{1cm}
&

 \begin{tikzpicture}[->,>=stealth',shorten >=1pt,auto,node distance=2.5cm,main node/.style={circle,draw,font=\sffamily\bfseries \small}, scale=1.35]
    % Draw axes
    \draw [->,thick] (0,-1.2)--(0,1.3) node (yaxis) [above] {$\widetilde{Q}(V_X(n_4))$}
        |- (-1.75,0)--(1.75,0) ;

   \node (a) at (0.524,0) {};
    \node (b) at (-1.007,0.513){};
    \node (c) at (-1.007,-0.513){} ;
    \node (d) at (1.490,0){} ;

\node (a1) at (0.524,0.3) {$F_1$};
    \node (b1) at (-1.007,0.813) {$E_1$};
    \node (c1) at (-1.057,-0.813) {$E_2$} ;
    \node (d1) at (1.490,0.3) {$E_0$};

   \fill[black] (a) circle (2.5pt);
\fill[black] (b) circle (2.5pt);
\fill[black] (c) circle (2.5pt);
\fill[black] (d) circle (2.5pt);

\path[every node/.style={font=\sffamily\small}]
(a) edge [dashed,->,bend right]  (d)
     (d) edge [bend right] (a)
  (c) edge [dashed,->,bend right] (b)
(b) edge [bend right] (c);
\end{tikzpicture}
\end{array}
$$

\bigskip

\hspace{-0.6cm} (iv) For $\underline{X=Bl_2(\mathbb{ P}^2)}$ the monodromy quivers are:

$$ %
$$
\bigskip

\hspace{-0.6cm} (v) For $\underline{ X=Bl_3(\mathbb{ P}^2)}$ the monodromy quivers are:

$$ %
$$

\bigskip

\hspace{-0.6cm} 5.2.2 \bf The Fano $\mathbb{P}^1$-bundles over $\mathbb{P}^2$ case \rm

\bigskip

\hspace{-0.6cm} (i) For $\underline{ X=\mathbb{P}(\mathcal{O}_{\mathbb{P}^2} \oplus \mathcal{O}_{\mathbb{P}^2}(1))}$ the monodromy quivers are:
$$ %
$$

\bigskip

\hspace{-0.6cm} (ii) For $\underline{ X=\mathbb{P}(\mathcal{O}_{\mathbb{P}^2} \oplus \mathcal{O}_{\mathbb{P}^2}(2))}$ the monodromy quivers are:

$$ %
$$

\hspace{-0.6cm} 5.2.3 \bf  The Fano $\mathbb{P}^2$-bundles over $\mathbb{P}^1$ case \rm

\bigskip

\hspace{-0.6cm} (i) For $\underline{ X=\mathbb{ P}(\mathcal{ O}_{\mathbb{P}^1} \oplus \mathcal{ O}_{\mathbb{P}^1} \oplus \mathcal{ O}_{\mathbb{P}^1} (1))}$ the monodromy quivers are:
$$ %
 $$

\hspace{-0.6cm} 5.2.4 \bf  The Fano $\mathbb{P}^1$-bundles over $\mathbb{P}^1 \times \mathbb{P}^1$ case \rm

\bigskip

\hspace{-0.6cm} (i) For $\underline{ X=\mathbb{P}(\mathcal{O}_{\mathbb{P}^1 \times \mathbb{P}^1} \oplus \mathcal{O}_{\mathbb{P}^1 \times \mathbb{P}^1}(1,1))}$ the monodromy quivers are:
$$ %
$$

\bigskip

\hspace{-0.6cm} (ii) For $\underline{ X=\mathbb{P}(\mathcal{O}_{\mathbb{P}^1 \times \mathbb{P}^1} \oplus \mathcal{O}_{\mathbb{P}^1 \times \mathbb{P}^1}(1,-1))}$ the monodromy quivers are:

$$ %
$$

\section{Discussion and Concluding Remarks}
\label{s:DCR}

\hspace{-0.6cm} In this work we showed examples of toric Fano manifolds $X$ which exhibit non-trivial relations between their small quantum cohomology $QH(X)$ and properties
of their derived category of coherent sheaves $\mathcal{D}^b(X)$. Concretely, between the Landau-Ginzburg solution scheme $Crit(f_X)$ (and its monodromies) and full strongly exceptional collections of
line bundles $\mathcal{E} \subset Pic(X)$ (and their quivers). The question is, of course, to which extent do these relations generalize to further semi-simple toric Fano manifolds?

\hspace{-0.6cm} Recall that the pair $((\mathbb{C}^{\ast})^n,f_X)$, where $f_X$ is the Landau-Ginzburg potential of $X$, is typically considered as the Hori-Vafa mirror of the toric Fano manifold $X$, see \cite{HV}. Kontsevich homological
mirror symmetry conjecture, in this setting, suggests that there is an equivalence of categories of the form $$ \mathcal{D}^b(X) \simeq \mathcal{D}^b(FS((\mathbb{C}^{\ast})^n,f_X)) $$ where $FS((\mathbb{C}^{\ast})^n, f_X)$ is the Fukaya-Seidel
category of $((\mathbb{C}^{\ast})^n, f_X)$, see \cite{Ko2}. We would like to conclude this work by mentioning a few remarks on relations to the framework of mirror symmetry:
\bigskip

\hspace{-0.6cm} \bf Remark 6.1 \rm (Seidel's vanishing cycles): In \cite{S} Seidel associates to an exact Morse fibration $f : E \rightarrow \mathbb{C}$ the derived
category of Lagrangian vanishing cycles $ \mathcal{D}^b(Lag_{vc}(f_X))$.  The definition of $ \mathcal{D}^b(Lag_{vc}(f_X))$ involves the choice of a basis of vanishing cycles $ \left \{ L_1,...,L_N \right \}$ in a non-singular fiber, which serves as a set of generators for $ \mathcal{D}^b(Lag_{vc}(f))$. Seidel shows the following equivalence of categories $\mathcal{D}^b(X) \simeq  \mathcal{D}^b(Lag_{vc}(f_X))$, in the case $X=\mathbb{P}^2$, where $f_X$ is the Landau-Ginzburg potential of $X$. The equivalence is obtained
via a specific choice of generators $\left \{ \widetilde{L}_0,\widetilde{L}_1,\widetilde{L}_2 \right \}$. Under the equivalence $ i : \mathcal{D}^b(Lag_{vc}(f_X)) \rightarrow \mathcal{D}^b(X)$ the collection $\left \{ i(\widetilde{L}_0),i(\widetilde{L}_1),i(\widetilde{L}_2) \right \}$ is
a full strongly exceptional collection. This method was later extended for further manifolds, for instance by Auroux, Katzarakov and Orlov for Del-Pezzo manifolds, see \cite{AKO}. However, in general, the resulting collections are not collections of \emph{line bundles}.

\bigskip

\hspace{-0.6cm} \bf Remark 6.2 \rm (The coherent-constructible correspondence): One of the developments in mirror symmetry for toric manifolds in recent years is the coherent-constructible
correspondence of Fang, Liu, Treumann and Zaslow, see \cite{FLTZ}. To a toric manifold $X$ the authors associate a Lagrangian submanifold $\Lambda_X \subset  (M_{\mathbb{R}}/M) \times N_{\mathbb{R}}$, defined in terms of
the toric data. One of the main geometric ingredients is the establishment of a relation between coherent sheaves on $X$ and constructible sheaves on $T^{\ast}T^{\vee}_{\mathbb{R}}:=(M_{\mathbb{R}}/M) \times N_{\mathbb{R}}$ with support
in $\Lambda_X$. Constructible sheaves, in turn, are related to the elements of a corresponding Fukaya category via a process of "microlocalization" due to works of Nadler-Zaslow, see \cite{NZ} and Nadler, see \cite{N}.

\hspace{-0.6cm} Denote by $\mathcal{P}rev(T^{\ast}X ; \Lambda)$ the category of perverse sheaves on $T^{\ast}X$ with support in a conical Lagrangian subvariety $\Lambda \subset T^{\ast}X$ (which is, by definition, a sub-category of the corresponding
category of constructible sheaves). A seminal result on the structure of $\mathcal{P}rev(T^{\ast}X ; \Lambda)$ due to S. Gelfand, MacPherson and Vilonen shows an equivalence of categories between $\mathcal{P}rev(T^{\ast}X ; \Lambda)$ and
 the category of representations of a quiver $\widetilde{Q}^{prev}(T^{\ast}X, \Lambda)$ which is, in turn, defined via monodromies, see \cite{GMV}. It is thus interesting  to ask whether the constructible dg-category $Sh_c(T^{\ast}T^{\vee}_{\mathbb{R}} ; \Lambda_X)$ studied by Fang, Liu, Treumann and Zaslow admits an analog quiver description. In view of the above, we would cautiously suggest that the analog monodromy map involved in the definition of such a quiver is the map $ M: \pi_1(L(\Delta^{\circ}) \setminus R_X; f_X) \rightarrow Aut(Crit(f_X))$.

\bigskip

\hspace{-0.6cm} \bf Acknowledgements: \rm Most of the results of this paper have been obtained during the author's postdoc at Neuch$\hat{\textrm{a}}$tel University. The author would 
like to thank Felix Schlenk for his ongoing support and scientific encouragement. The author would also like to express his gratitude to Paul Biran for 
 ongoing support and for introducing him to the subject of quantum cohomology.

\end{document}